%% file: chromatic.tex
\documentclass{amsart}

\usepackage{amsmath}
\usepackage{amsthm}
\usepackage{amsfonts}
\usepackage{amssymb,latexsym}
\usepackage[all]{xy}
\usepackage{graphicx}
\usepackage[mathscr]{eucal}
\usepackage{verbatim}

\theoremstyle{plain}
    \newtheorem{thm}{Theorem}[section]
    \newtheorem{prop}[thm]{Proposition}
    \newtheorem{lemma}[thm]{Lemma}
    \newtheorem{conj}[thm]{Conjecture}
    \newtheorem{cor}[thm]{Corollary}

\theoremstyle{definition}
    \newtheorem{defn}[thm]{Definition}

\theoremstyle{remark}
    \newtheorem{rem}[thm]{Remark}
    \newtheorem{example}[thm]{Example}

\newtheorem*{namedtheorem}{\theoremname}
\newcommand{\theoremname}{testing}

\newcommand{\A}{\mathcal A}
\newcommand{\B}{\mathcal B}
\newcommand{\C}{\mathcal C}
\newcommand{\D}{\mathcal D}

\newcommand{\F}{\mathcal F}

\newcommand{\T}{\mathcal T}

\newcommand{\stk}[1]{\stackrel{#1}{\longrightarrow}}
\newcommand{\la}{\ensuremath{\langle}}
\newcommand{\ra}{\ensuremath{\rangle}}
\newcommand{\ints}{\ensuremath{\mathbb{Z}}}

\newcommand{\rats}{\ensuremath{\mathbb{Q}}}

\newcommand{\ftwo}{\mathbb{F}_2}
\newcommand{\rar}{\ensuremath{\rightarrow}}

\newcommand{\n}{\noindent}
\newcommand{\lspace}{\hspace{2pt}}

\DeclareMathOperator{\Supp}{Supp}
\DeclareMathOperator{\Spec}{Spec}
\DeclareMathOperator{\Hom}{Hom}

\DeclareMathOperator{\proj}{proj}

\DeclareMathOperator{\Ker}{Ker}

\DeclareMathOperator{\plainlim}{lim}
\newcommand{\clim}[1]{\ensuremath{\underset{#1}{\plainlim}\lspace}}

\DeclareMathOperator{\plaincolim}{colim}
\newcommand{\colim}[1]{\ensuremath{\underset{#1}{\plaincolim}\lspace}}

\newcommand{\rlim}[1]{\ensuremath{\underset{#1}{\text{Rlim}}\lspace}}

\DeclareMathOperator{\plainhocolim}{hocolim}
\newcommand{\hocolim}[1]{\ensuremath{\underset{\scriptscriptstyle #1}{\plainhocolim}\lspace}}

\newcommand{\sus}[1]{\ensuremath{\Sigma ^{#1}}}

\newcommand{\bp}[1]{\ensuremath{ \text{BP} \langle {#1} \rangle}}
\newcommand{\BP}{\ensuremath{ \text{BP}}}

\DeclareMathOperator{\Pic}{Pic}
\DeclareMathOperator{\im}{Im}

\begin{document} 

\title{Refining thick subcategory theorems}

\author{Sunil K. Chebolu} 
\address {Department of Mathematics,
          University of Western Ontario,
          London, Ontario, Canada, N6A 5B7.}
\email{schebolu@uwo.ca}


\begin{abstract} We use a $K$-theory recipe of Thomason to obtain classifications of triangulated subcategories via
refining some standard thick subcategory theorems. We apply this recipe to the full subcategories of finite objects in the  
derived categories of rings and the stable homotopy category of spectra. This gives, in the derived categories, a complete classification
of the triangulated subcategories of perfect complexes over some commutative rings. In the stable homotopy category of spectra we obtain only a partial classification of
the triangulated subcategories of the finite $p$-local spectra. We use this partial classification to study the lattice of triangulated subcategories.
This study gives some  new evidence to a conjecture of Adams that the thick subcategory $\C_2$ can be generated by iterated cofiberings of the Smith-Toda complex.
We also discuss several consequences of these classification theorems.
\end{abstract}

\maketitle

\section{Introduction}

Classifying various subcategories of triangulated categories like the derived categories and the homotopy category of spectra has been an active area and has
proved to be  extremely useful in the study of global problems in (stable) homotopy theory.
Several mathematicians brought to light many amazing a priori different theories by classifying various subcategories of triangulated categories.
Following the seminal work of Devinatz, Hopkins, and Smith \cite{dhs} in stable homotopy theory, this line of research was initiated 
by Hopkins in the 80s. In his famous 1987 paper \cite{Ho}, Hopkins classified the thick subcategories (triangulated subcategories that are closed under retractions) 
of the finite $p$-local spectra and those of perfect complexes over a noetherian ring. He showed that thick subcategories of the finite spectra are determined by the 
Morava $K$-theories and those of perfect complexes by the prime spectrum of the ring. These results have had tremendous impacts in their respective fields.
The thick subcategory theorem for finite spectra played a vital role in the study of nilpotence and periodicity.  For example, using this theorem
Hopkins and Smith \cite{hs} were able to settle the class-invariance conjecture of Ravenel \cite{rav} which classified the Bousfield classes of finite spectra. 
Similarly the thick subcategory theorem for the derived category  establishes a surprising connection between stable homotopy theory and algebraic 
geometry; using this theorem one is able to recover the spectrum of a ring from the homotopy structure of its derived category!  
These ideas were later pushed further 
into the world of derived categories of rings and schemes by Neeman \cite{Ne} and Thomason \cite{Th}, and into modular representation theory by Benson, Carlson and 
Rickard \cite{bcr}. Motivated by the work of Hopkins \cite{Ho}, Neeman \cite{Ne} classified the Bousfield classes and localising subcategories in the derived category of a 
noetherian ring. In modular representation theory, the Benson-Carlson-Rickard classification of the thick subcategories of stable modules over group algebras has led to 
some deep structural information on the representation theory of finite groups. Finally the birth of axiomatic stable homotopy theory \cite{mps} in the mid 90s encompassed 
all these various theories and ideas and studied them all in a more general framework. With all these developments over the last 30 years, the importance of triangulated 
categories in modern mathematics is by now abundantly clear.

The goal of this paper is to classify the triangulated subcategories analogous to the aforementioned classifications of thick subcategories. To motivate this project further,
let us consider the following question. Let $\T$ be a triangulated category and $X$ and $Y$ be two objects in $\T$.

\vskip 3mm \n
\textbf{Question A.} When can $Y$ be generated from $X$ using \emph{cofibrations} and \emph{retractions}? \\

Knowledge of the thick subcategories of $\T$ will help us  answer this question. For example,  if $\T$ is the category of perfect complexes over a noetherian ring or the 
category of finite $p$-local spectra, then we know (from the Hopkins thick subcategory theorems) that $Y$ can be generated from $X$ using cofibrations and retractions if
and only if $\Supp(Y) \subseteq \Supp(X)$. (When $X$ is a perfect complex, $\Supp(X)$ is the set of primes $p$ for which $X \otimes R_p \ne 0$; and when 
$X$ is a finite $p$-local spectrum, by  $\Supp(X)$ we mean the chromatic support of $X$, i.e., the set of non-negative integers $n$ for which $K(n)_*X \ne 0$.)
This is quite remarkable because often ``support''  is a computable invariant while cofibrations and retractions can be extremely hard to understand. Now let us ask a
even more subtle and stringent question.

\vskip 3mm \n
\textbf{Question B.} When can $Y$ be generated from $X$ using  \emph{cofibrations} alone? \\

We say that this is a stringent question because it is much harder, in general, to work with cofibrations alone. For example, take $X=M(p)$ and $Y=M(p^2)$, clearly
$\Supp(M(p)) = \Supp(M(p^2))$, therefore  $M(p)$ can be generated from $M(p^2)$ using cofibrations and retractions. However, it is impossible, 
as we will see, to generate $M(p)$ from $M(p^2)$ just using cofibrations.  On the other hand, $M(p^2)$ can be generated from $M(p)$ using cofibrations: There is a cofibre 
sequence
\[ \Sigma^{-1} M(p) \rar M(p) \rar M(p^2).\]
So having a classification of triangulated subcategories will help us answer Question B. This is our motivation for classifying triangulated subcategories.
In order to classify the triangulated subcategories, we use a $K$-theoretic approach of Thomason to refine Hopkins's
thick subcategory theorems, both for finite spectra and perfect complexes. Note that, in both these categories, the support condition ($\Supp(Y) \subseteq \Supp(X)$) is 
a necessary condition for the generation in  Question B.  Thomason's $K$-theory technique helps us to come up with a sufficient condition! More precisely, we construct
universal Euler characteristic functions $\chi$ on $\T$, and show that $Y$ can be generated from $X$ using cofibrations if and only if 
 $\Supp(Y) \subseteq \Supp(X)$ and  $\chi(X)$ divides $\chi(Y)$; see corollary \ref{cor:testforp-localspectra}, corollary \ref{cor:testforp-torsionspectra}, 
corollary \ref{cor:questionBPID}, and corollary \ref{cor:questionBArtin}.

We also discuss the following conjecture of Adams. In his last (unpublished) paper \cite{MayTho},   
Adams conjectured that the Smith-Toda complex (which is known to exist 
at odd primes) generates the thick subcategory $\C_2$ by cofibrations. (Note that this is clearly possible if we allow retractions.) Using $\text{BP}$-based homology theories,
we construct Euler characteristic functions $\chi_n$ on the thick subcategory $\C_n$:
\[ \chi_n(X) := \sum_i (-1)^i \log_p |\bp{n-1}_i X| .\]
These Euler characteristic functions give us triangulated subcategories $\C_n^k$ of $\C_n$: 
   \[ \C_n^k = \{ X \in \C_n \; | \; \chi_n(X) \equiv 0 \ \text{mod}  \, l_n k \},\]
where $l_n$ is the smallest positive value of $\chi_n$ on $\C_n$.
We then show (using a Bockstein spectral sequence calculation) that, for all $n$ and $k$, 
\[ \C_{n+1} \subsetneq  \C_n^k\subsetneq \C_n.\]
This gives only some evidence to Adams's conjecture (when $n=1$ the above inclusions follow trivially from Adams conjecture). The conjecture, however, remains open.

This paper is organised as follows. We set up our categorical stage in the next section. We begin with a quick recap of Grothendieck groups and some categorical 
definitions. We then explain a $K$-theoretic technique of Thomason which will be used in the later sections to obtain classifications of triangulated subcategories. 
We also revise Thomason's theorem for triangulated categories that are equipped with a nice product.
We apply these techniques to spectra (section 3) and derived categories (section 4).  As consequences of these classifications, we will study various
structural results on triangulated subcategories. We also record some questions that come up along the way.

\vskip 4mm \n
\textbf{Acknowledgements} I would like to thank my advisor John Palmieri for the innumerable conversations that I had with him on this material. Among others, I thank
Mark Hovey, for bringing to my attention the work of Adams on Universal Euler characteristics; Steve Mitchell and Hal Sadofsky, for some interesting discussions on the 
subject. I am also grateful to an anonymous referee for giving me many valuable suggestions.


\section{Thomason's $K$-theory recipe}

\subsection{Grothendieck groups}
We begin by recalling some definitions and results from
\cite{Th}. Let $\T$ denote a  triangulated category that is
essentially small (i.e., it has only a set of isomorphism classes of
objects). Then the \emph{Grothendieck group} $ K_0(\T) $ is defined to
be the free abelian group on the isomorphism classes of $\T$ modulo
the Euler relations $[B]=[A]+[C]$, whenever $A\rightarrow B
\rightarrow C \rightarrow \Sigma A $ is an exact triangle in $\T$
(here $[X]$ denotes the element in the Grothendieck group  that is
represented by the isomorphism class of the object $X$). This is
clearly an abelian group with $[0]$ as the identity element and
$[\Sigma X]$ as the inverse of $[X]$. The identity
$[A]+[B]=[A \amalg B]$ holds in the Grothendieck group. Also note that any
element of $K_0(\T)$ is of the form $[X]$ for some $X\in \T$. All
these facts follow easily from the axioms for a triangulated category. 

Grothendieck groups have the following universal property: Any map $\sigma$ from the set of  isomorphism
classes of $\T$  to an abelian group $G$ such that the Euler relations hold in $G$
factors through a unique  homomorphism $ f: K_0(\T)
\rightarrow G $; see the diagram below.
\[
\xymatrix{ {\{\mbox{isomorphism classes of objects in}\, \T\}} \ar[d]^\pi \ar[dr]^
\sigma \\  {K_0(\T)} \ar@{..>}[r]^f & G }
\]
Therefore the map $\pi$ will be called the \emph{universal Euler characteristic function}.
$K_0(-)$ is clearly a covariant functor from the category of small triangulated categories to the 
category of abelian groups.

Unless stated otherwise, all subcategories in this paper are assumed to be \emph{full} subcategories.

\begin{defn} \label{defn:small}
An object $X$ in a triangulated category $\T$ is \emph{small} if the natural map 
\[ \underset{\alpha \in \Lambda}{\bigoplus}\, \Hom(X, A_{\alpha}) \rightarrow  \Hom(X, \,\underset{\alpha \in \Lambda}{\coprod} A_{\alpha})\]
is an isomorphism for all set indexed collections of objects $A_{\alpha}$ in $\T$. 
(Some authors call such objects as finite or compact objects.)
\end{defn}

\begin{defn} \label{defn:thick}
A triangulated subcategory $\C$ of $\T$ is said to be a \emph{thick subcategory} if it is closed under
retractions, i.e., given a commuting diagram
\[
\xymatrix{
 B \ar[r]^i \ar@/_20pt/[rr]_= & A \ar[r]^r & B
}
\]
in $\T$ such that $A$ is an object of $\C$, then so is $B$. Since retractions split in a triangulated category, 
this property of $\C$ is equivalent to saying that $\C$ is closed under 
direct summands: $A \amalg B \in \C  \Rightarrow  A \in \C$ and 
$B \in \C$. 
\end{defn}

\begin{example} The full subcategory of small objects in any triangulated category is thick. 
\end{example}

\begin{defn} \label{defn:dense}
We say that a triangulated subcategory $\C$ is \emph{dense} in $\T$ if every object in $\T$ is a direct 
summand of some object in $\C$. 
\end{defn}

The following theorem due to Thomason \cite{Th} is the foundational theorem
that motivated this paper.

\begin{thm} \cite[Theorem 2.1]{Th} \label{main} Let $\T$ be an essentially small triangulated category. 
Then there is a natural order preserving bijection between the posets,

\begin{center}
\{dense triangulated subcategories $\A$ of $\T$ \}
$\overset{f} {\underset{g}{\rightleftarrows}}$  \{subgroups $H$ of
$K_0(\T)$ \}.
\end{center}

\noindent The map $f$ sends $\A$ to image of the map $K_0(\A) \rar K_0(\T)$, 
and the map $g$ sends $H$ to the full subcategory of all objects $X$
in $\T$ such that $[X] \in H$.
\end{thm}

\subsection{Thomason's $K$-theory recipe} \label{ss:recipe}

The importance of Thomason's theorem \ref{main} can be seen from the following simple observation.
\emph{Every triangulated subcategory $\mathcal{A}$ of $\T$ is dense in a unique thick
subcategory of $\T$} --  the one obtained by taking the intersection of all the thick subcategories 
of  $\T$ that contain $\mathcal{A}$. This observation in conjunction with theorem \ref{main} gives the following brilliant 
recipe of Thomason to the problem of classifying the triangulated subcategories of $\T$:

\begin{enumerate}
\item{Classify the thick subcategories of $\T$.}
\item{Compute the Grothendieck groups of all thick subcategories.}
\item{Apply Thomason's theorem \ref{main} to each thick subcategory of $\T$.}
\end{enumerate}

We will apply this recipe to the categories of small objects in some stable homotopy categories like  
the stable homotopy category of spectra and the derived categories of rings. 
The stable homotopy structure on these categories will guide us while applying this recipe and 
consequently we will derive some structural information on these categories.

\subsection{Grothendieck Ring}
\noindent Through out this subsection $\T$ will denote a tensor triangulated 
category that is essentially small. Now making use of the available smash
product, we want to define a ring structure on the Grothendieck group. This can be done in a very natural and obvious way.

\begin{defn} If $[A]$ and $[B]$
are any two elements of $K_0(\T)$, then define $[A] [B] := [A \wedge B] .$
\end{defn}

This can be easily shown to be a well-defined operation and endows $K_0(\T)$ with the structure of a commutative ring. The Grothendieck class 
of the unit object $[S]$ serves as the identity element in the ring.

The following three definitions are motivated by their analogues in commutative ring theory.

\begin{defn} A full triangulated
subcategory $\C$ of $\T$ is said to be \emph{$\otimes$-closed} or a
\emph{triangulated ideal} if for all $A \in \T$ and for all $B \in \C $, $B
\wedge A \in \C$. A full triangulated ideal is said to be respectively
thick or dense if it is such as a triangulated subcategory.
\end{defn}

\begin{defn} A full triangulated ideal $\B$ of $\T$ is said to be \emph{prime} if
for all $X$, $Y \in \T$ such that  $X \wedge Y \in \B$,  either $X
\in \B$ or $ Y \in \B$. Similarly $\B$ is said to be  a \emph{maximal} in $\T$ if there
is no triangulated ideal $\A$ such that $\B \subsetneq \A \subsetneq \T$.
\end{defn}

\begin{defn} We say that a triangulated category $\A$ is a \emph{triangulated module} over a tensor triangulated category $\T$
if there is a triangulated functor
\[ \phi: \T \times \A \rar \A\]
that is covariant and exact in each variable and satisfies the obvious unital and associative conditions.
A full triangulated subcategory $\B$ of $\A$ is a \emph{triangulated submodule} of $\A$ if the functor
$\phi$ maps $\T \times \B \rar \B$. Note that in this situation, $K_0(A)$ becomes a $K_0(\T)$-module, and $K_0(B)$  a $K_0(\T)$-submodule.
\end{defn}
 
We need the following lemma to upgrade Thomason's theorem to tensor triangulated categories.

\begin{lemma} \cite[Lemma 2.2]{Th} Let $\A$ be a dense triangulated subcategory of an essentially small triangulated category
$\T$. Then for any object $X$ in $\T$, one has that $X \in \A$ if and
only if $[X] = 0 \in K_0(\T)/ \im (K_0(\A) \rar K_0(\T)).$
\end{lemma}

The following results are now expected.

\begin{thm} \label{imain}
Let $\T$ be a tensor triangulated category that is essentially
small. Then, under Thomason's bijection (theorem \ref{main})
\begin{center}
\{dense triangulated subcategories of $\T$ \} $\longleftrightarrow$ \{subgroups of $K_0(\T)$\}, 
\end{center}
we have the following correspondence.
\begin{enumerate} 
\item The dense triangulated ideals correspond precisely to the ideals of
the ring $K_0(\T)$.
\item The dense prime triangulated ideals correspond precisely to the
prime ideals of $K_0(\T)$.
\item The dense maximal triangulated ideals correspond precisely to
the maximal ideals of $K_0(\T)$.
\end{enumerate}
\end{thm}
\begin{proof} Except possibly the second statement, everything else is straightforward.

\noindent 1. Let $I$ be any ideal in $K_0(\T)$. Then the
corresponding dense triangulated subcategory is  $\T_{I} = \{X \in \T :[X] \in I \}$.  
Now for $A \in \T$ and $B \in A_I$, note that $[B \wedge A]=[B][A] \in I$ ($I$ being an ideal) and
therefore $B \wedge A \in \T_{I}$. This shows that $\T_{I}$ is a triangulated ideal.
The other direction is equally easy.

\noindent 2. Suppose $H$ is a prime ideal in $K_0(\T)$.  The
corresponding dense triangulated ideal is given by $\B = \{ X : [X] \in H \}$. 
Now if $A \wedge B \in  \B$, then by definition of $\B$, we have
$[A\wedge B] \in H $, or equivalently $[A][B] \in H $. Now primality of
$H$ implies that either $[A] \in H$ or $[B] \in H$, which means 
either $A \in \B $ or $ B \in \B$. For the other direction, suppose
$\B$ is dense prime triangulated ideal of $\T$. The
corresponding subgroup $H$ is the image of the map $K_0(\B) \rar K_0(\T) $. 
Now suppose the product of two elements $[A]$ and $[B]$ belongs to $H$. Then we have $[A \wedge B] \in H$ or
equivalently  $[A \wedge B] = 0 $ in $K_0(\T)/ H$.  By the above lemma, we then have
 $A \wedge B \in \B$.  Since $\B$ is prime,
this implies that either $A \in \B$ or $B \in \B$, or equivalently $[A] \in H$ or $[B] \in H$.

\noindent 3. This follows directly from the fact that Thomason's bijection is
a map of posets. 
\end{proof}

In the same spirit we get the following result for triangulated modules over tensor triangulated categories. We leave the
proof which is similar to the above theorem as an easy exercise to the reader.

\begin{thm} \label{refree-main}
Let $\A$ be a triangulated module over a tensor triangulated category $\T$.
Then, under Thomason's bijection (theorem \ref{main})
\begin{center}
\{dense triangulated subcategories of $\A$ \} $\longleftrightarrow$ \{subgroups of $K_0(\A)$\}, 
\end{center}
we have the following correspondence.
\begin{enumerate} 
\item The dense triangulated submodules of $\A$ correspond precisely to the $K_0(\T)$-submodules of
$K_0(\A)$.
\item The dense maximal triangulated submodules correspond precisely to
the maximal $K_0(\T)$-submodules of $K_0(\A)$.
\end{enumerate}
\end{thm}

\noindent
\section{Triangulated subcategories of finite spectra} \label{se:K_0groups}

In this section we apply Thomason's theorem \ref{main} to the category $\F_p$ of finite $p$-local spectra.
So we begin by recalling the celebrated thick subcategory theorem of Hopkins and Smith which is a key ingredient in Thomason's recipe for classifying the
triangulated subcategories of $\F_p$. We then examine the Grothendieck groups of these thick subcategories. 
Our partial knowledge of these groups gives us a family of triangulated subcategories of
$\F_p$. We then study the lattice of these subcategories in subsection \ref{se:lattice} using some spectral sequences. This study gives new evidence to a
conjecture of Adams. We end the section with some questions.

For each non-negative integer $n$, there is a field 
spectrum called the $n$th Morava $K$-theory $K(n)$, whose coefficient ring is
$\mathbb{F}_p[v_n,v_n^{-1}]$ with $|v_n| = 2(p^n-1)$. According to the thick subcategory theorem of Hopkins-Smith 
these Morava $K$-theories determine the thick subcategories of $\F_p$. More precisely:

\begin{thm} \cite{hs} For each positive integer $n$,
let $\C_n$ denote the full subcategory of all finite $p$-local
spectra  that are $K(n-1)$-acyclic. Then a non-zero subcategory $\C$ of $\F_p$ is 
thick if and only if $\C = \C_n$ for some $n$. Further these thick subcategories
give a nested decreasing filtration of $\F_p$ \cite{rav, mit}: 
\[ \cdots \C_{n+1} \subsetneq \C_n \subsetneq \C_{n-1} \cdots \subsetneq \C_1 \subsetneq \C_0 (= \F_p).\]
\end{thm}

A property $P$ of finite spectra is said to be \emph{generic} if the collection of all spectra in $\F_p$ which satisfy the given property $P$ is $\C_n$ for some $n$.
A spectrum is said to be of type-$n$ if it  belongs to $\C_n - \C_{n+1}$. For example, the sphere spectrum $S$ is of type-$0$, the mod-$p$ Moore spectrum
$M(p)$ is of type-$1$, etc.

The problem of computing the Grothendieck groups of thick subcategories of the
finite $p$-local spectra was first considered, to our knowledge, by Frank Adams.
This appeared in an unpublished manuscript \cite[Page 528-529]{MayTho} of Adams on the work of Hopkins.
We begin with a recapitulation of Adams's work and then bring it to this new context of classifying
triangulated subcategories. We should also point out that our Euler characteristic functions are simpler
than the ones considered by Adams.

\subsection{$\C_0$ - Finite $p$-local spectra} 
We begin with the fundamental notion of Euler characteristic of a finite spectrum. Recall that if $X$ is any finite $p$-local spectrum, 
the Euler characteristic of $X$ with rational coefficients  is given by  
\begin{equation} \label{E:chi_0}
\chi_0 (X) = \sum_{i= -\infty}^{\infty} (-1)^i \dim_{\rats} \;H\rats_i(X).
\end{equation}

Since $X$ is a finite spectrum it has homology concentrated only in some finite range and therefore this is a
well-defined function.

\begin{example} For every non-negative integer $m$, define a full subcategory $\C_0^m$ of $\C_0 \, (= \F_p)$ as:
 \[\C_0^m =\{ X \in \C_0 : \chi_0 (X)\equiv 0 \, \text{mod} \,m \} .\] 
It is an easy exercise to verify that these are all dense triangulated subcategories of $\C_0$.
\end{example}

\begin{prop} 
A triangulated subcategory $\C$ of $\C_0 (= \F_p)$ is dense in  $\C_0$ if and only if $\C = \C_0^m$ for some $m$.
\end{prop}

\begin{proof} Let $S$ denote the $p$-local sphere spectrum and recall that $\C_0$ can be generated by iterated cofiberings of the sphere
spectrum $S^0$. Now using the Euler relations, it is clear that the 
Grothendieck group of $\C_0$ is a cyclic group generated by $[S]$. 
The map                                    
\[\chi_0 :\{\mbox{isomorphism classes of }  \C_0 \}  \rightarrow \mathbb{Z}\] 
which sends an isomorphism class to its Euler characteristic shows that $K_0(\C_0)$ is isomorphic to $\ints$. 
The proposition now follows by invoking theorem \ref{main}.  
\end{proof}

\begin{prop} $K_0(\C_0) \cong \mathbb{Z}$ as commutative rings.
\end{prop}
\begin{proof} First note that the smash product on $\C_0$ induces a ring structure on $K_0(\C_0)$. 
 Now the isomorphism of abelian groups $\phi_0: K_0(\C_0) \rar \mathbb{Z}$ (from the proof of the previous proposition)
maps the Grothendieck class of the wedge of $n$-copies of the sphere spectrum  to $n$.  Since the smash product distributes over the wedge, we 
conclude that $\phi_0$ is a ring isomorphism.  
\end{proof}

We now give some nice consequences of this proposition.

\begin{cor} Every triangulated subcategory of $\C_0$ is a triangulated ideal.
\end{cor}

\begin{proof} Let $\C$ be a triangulated subcategory of $\C_0$, and let $\C_n$ denote the unique thick subcategory in which
$\C$ is dense.  Now note that $\C_n$ is a triangulated module over $\C_0$, therefore $K_0(\C_n)$ is a $K_0(\C_0)$-module. 
We have seen that $K_0(\C_0) \cong \mathbb{Z}$ as rings. So we apply theorem \ref{refree-main}
and observe the simple fact that subgroups of the abelian group $K_0(\C_n)$ are precisely the $\ints$-submodules of $K_0(\C_n)$. 
\end{proof}

\begin{cor} If $\B$ is a dense triangulated subcategory of $\F_p$, then $\B$ is prime if and only if $\B =\C_0^0 $ or $\C_0^p$ for some
prime number $p$.  In particular, a dense triangulated subcategory of $\F_p$ is maximal if and only if it is $\C_0^p$ for some prime $p$.
\end{cor}

\begin{proof}  By theorem \ref{imain}, we have a correspondence between the prime (maximal) dense triangulated subcategories
of $\F_p$ and the prime (maximal) ideals of $K_0(\F_p) \,(\cong \ints)$. The corollary now follows by comparing 
the prime ideals and maximal ideals of $\mathbb{Z}$.
\end{proof}

\begin{cor} \label{cor:testforp-localspectra} Let $\Delta(X_1,X_2,\cdots,X_k)$ denote the triangulated subcategory generated by 
spectra $X_1, X_2, \cdots,$ and $X_k$. Then we have the following.
\begin{enumerate}
\item If $A$ is a  type-$0$ spectrum, then $\Delta(A)$ consists of all spectra $X$ in $\C_0$ for which $\chi_0(X)$ is divisible by $\chi_0(A)$.
\item If $A$ is a type-$0$ spectrum and $B$ is a spectrum in $\C_0$, then $B$ can be generated from $A$ using cofibre sequences if and only if $\chi_0(B)$  is divisible by 
$\chi_0(A)$.
\item If $A$ is a type-$0$ spectra and $B$ is a spectrum in $\C_0$, then $\Delta(A,B)$ consists of all spectra $X$ in $\C_0$ for which 
$\chi_0(X)$ is divisible by the 
highest common factor of $\chi_0(A)$ and $\chi_0(B)$.
\end{enumerate}
\end{cor}

\begin{rem} The above corollaries  are some interesting structural results on the subcategories of $\F_p$. It is not clear how one would
establish such results without using this $K$-theory approach of Thomason.  
\end{rem}

\subsection{$\C_1$ - Finite $p$-torsion spectra.}
Having classified the dense triangulated subcategories of $\C_0$, we now
look at the thick subcategory $\C_1$ consisting of the finite $p$-torsion spectra.
A potential candidate which might generate $\C_2$ as a triangulated category is the Moore spectrum $M(p)$.
So before going any further, it is natural to ask if $M(p^2)$ can be generated from $M(p)$ by iterated cofiberings. The fact that this is possible can be seen
easily using the octahedral axiom which gives the following commutative diagram of cofibre sequences. 
\[
\xymatrix{
S \ar[r]^p  \ar@{=}[d]  & S \ar[r] \ar[d]^p & M(p) \ar@{.>}[d] \\
S \ar[r]^{p^2} \ar[d]^0 & S \ar[r] \ar[d] & M(p^2) \ar@{.>}[d] \\
0 \ar[r]& M(p)  \ar[r] &  M(p)
}
\]
The exact triangle in the far right is the desired cofibre sequence.
Inductively, it is easy to see that $M(p^i)$ can be generated from $M(p)$ using cofibre sequences.
This example motivates the next proposition.

\begin{prop} \label{prop:steve} Every spectrum in $\C_1$ can be generated by iterated cofiberings of $M(p)$.
\end{prop}

\begin{proof} First observe that the integral homology of any spectrum in $\C_1$ consists of 
finite abelian $p$-groups. Further these spectra have  homology concentrated only in a finite range, so we induct on $|H\mathbb{Z}_*(-)|$. 
If $|H\mathbb{Z}_*(X)| = 1 $, that means $X$ is a trivial spectrum and it is obviously generated by $M(p)$.  So assume that 
$|H\mathbb{Z}_*(X)| > 1$ and let $k$ be the smallest integer such that $H\mathbb{Z}_k(X)$ is non-zero. 
Then by the Hurewicz theorem, we know that $\pi_k(X) \cong H\mathbb{Z}_k(X)$. 
So pick an element of order $p$ in $H\ints_k (X)$ and represent it by a map 
$\alpha: S^k \rightarrow X$. Since $p\alpha = 0$, the composite $S^k \stk{p} S^k \stk{\alpha} X$ 
is zero and hence the map $\alpha$ factors through $\Sigma^k M(p)$.  This gives the following diagram 
where the vertical sequence is a cofibre sequence.
\[
\xymatrix{
\Sigma^k S \ar[drr]_0 \ar[r]^p & \Sigma^k S \ar[dr]^{\alpha} \ar[r]  & \Sigma^k M(p) \ar@{.>}[d]^{\alpha'}\\
& & X \ar[d] \\
& & Y 
}
\]
It is easily seen that $H\mathbb{Z}_i(\alpha')$ is non-zero (and hence injective) if $i=k$, and is zero otherwise.
Therefore, from the long exact sequence in integral homology induced by the vertical cofibre 
sequence, it follows that $|H\mathbb{Z}_*(Y)| = |H\mathbb{Z}_*(X)| - p$. The induction
hypothesis tells us  that $Y$ can be generated by cofibre sequences using $M(p)$ and then the above vertical 
cofibre sequence tells us that $X$ can also be generated by cofibre sequences using $M(p)$.
\end{proof}

Note that the regular Euler characteristic is not good for spectra in $\C_1$ because these
spectra are all rationally acyclic and therefore their Euler characteristic is always
zero (no matter which field coefficients we use). The integral homology of these spectra 
consists of finite abelian $p$-groups and that motivates the following definition.

\begin{defn} For any spectrum $X$ in $\C_1$, define 
\begin{equation} \label{E:chi_1}
\chi_1(X) := \sum_{i = -\infty}^{\infty} (-1)^i \log_p |H\mathbb{Z}_i X|,
\end{equation}
and for every non-negative integer $m$, define a full subcategory 
\[\C_1^m := \{ X \in \C_1 : \chi_1(X) \equiv 0\ \text{mod} \, m \}.\]
\end{defn}

\begin{thm} \label{thm:steve}A triangulated subcategory $\C$ of $\C_1$ is dense in $\C_1$ if and only if $\C = \C_1^m$ for some non-negative integer $m$.
\end{thm}

\begin{proof} The function $\log_p |-|$ is clearly an additive function on the abelian category of
finite abelian $p$-groups, i.e., whenever $0 \rightarrow A \rightarrow B \rightarrow C \rightarrow 0$
is a short exact sequence of finite abelian $p$-groups, we have $\log_p |B| = \log_p |A| + \log_p |C|$.
Using this it is elementary to show that if we have a bounded exact sequence of finite abelian $p$-groups,
then the alternating sum of $\log_p|-|$'s is zero. This applies, in particular, to the long exact sequence
in integral homology for a cofibre sequence in $\C_1$. In other words, $\chi_1(-)$ respects the
Euler relations in $\C_1$. So we get an induced a map $\phi_1 : K_0(\C_1) \rar \ints$. From  
proposition \ref{prop:steve}, we know that $K_0(\C_1)$ is a cyclic group generated by $[M(p)]$. 
Since $\phi_1([M(p)]) = 1$, it follows that $\phi_1$ is an isomorphism. The given classification
is now clear from theorem \ref{main}
\end{proof}

\begin{cor} A triangulated subcategory $\C$ of $\C_1$ is a maximal triangulated subcategory  if and only if $\C = \C_1^p$ for some prime number $p$.
\end{cor}

\begin{proof} Since $K_0(\C_1) \cong \ints$,  the corollary follows (as before) by observing that the maximal subgroups of $\ints$ are precisely the subgroups
generated by prime numbers.
\end{proof}

Now we state the analogue of corollary \ref{cor:testforp-localspectra} for finite $p$-torsion spectra.

\begin{cor} \label{cor:testforp-torsionspectra} Let $\Delta(X_1,X_2,\cdots,X_k)$ denote the triangulated subcategory generated by 
spectra $X_1, X_2, \cdots,$ and $X_k$. Then we have the following.
\begin{enumerate}
\item If $A$ is a  type-$1$ spectrum, then $\Delta(A)$ consists of all spectra $X$ in $\C_1$ for which $\chi_1(X)$ is divisible by $\chi_1(A)$.
\item If $A$ is a type-$1$ spectrum and $B$ is a spectrum in $\C_1$, then $B$ can be generated from $A$ using cofibre sequences if and only if $\chi_1(B)$  is divisible by 
$\chi_1(A)$.
\item If $A$ is a type-$1$ spectra and $B$ is a spectrum in $\C_1$, then $\Delta(A,B)$ consists of all spectra $X$ in $\C_1$ for which 
$\chi_1(X)$ is divisible by the 
highest common factor of $\chi_1(A)$ and $\chi_1(B)$.
\end{enumerate}
\end{cor}

We now examine the Grothendieck group of the Verdier Quotient category $\C_0/\C_n$. We begin with some generalities.
If $\A$ is any thick subcategory of a triangulated category $\C$, then we have natural functors fitting
into an exact sequence
\[ \A \rightarrow \C \rightarrow \C/\A \rightarrow 0,\]
where the first functor is the inclusion functor and the second one is the quotient functor into the
Verdier quotient. Applying $K_0(-)$ to the above sequence induces an exact sequence \cite[Page 355, Proposition 3.1]{sga5}, 
\begin{equation} \label{eq:K_0exactseq}
 K_0(\A) \rightarrow K_0(\C) \rightarrow K_0(\C/\A) \rightarrow 0. 
\end{equation}
Moreover, if $\C$ is a tensor triangulated category and $\A$ is a thick ideal in $\C$, then $K_0(\C/\A)$ is a $K_0(\C)$-algebra.
The first map in the above exact sequence is in general not injective. Here is an example where such a map
fails to be injective. If $\A = \C_n$ ($n \ge 1$) and $\C = \C_0$, then the inclusion functor 
$\C_n \hookrightarrow \C_0$  induces a map 
\[K_0(\C_n) \rightarrow K_0(\C_0) \cong \mathbb{Z}.\] 
This map is evidently the zero map because the Euler characteristic $\chi_0$ applied to any rational acyclic gives zero. 
Exactness of the sequence $K_0(\C_n) \rar K_0(\C_0) \rar K_0(\C_0/\C_1) \rar 0$ gives the following corollary.

\begin{cor} $K_0(\C_0/\C_n) \cong \mathbb{Z}$ as commutative rings, and  is generated by the Grothendieck class of the image of
the $p$-local sphere spectrum under the quotient functor $\C_0 \rightarrow \C_0/\C_n$.
\end{cor}

Since the thick subcategories of $\F_p$ are all nested ($\C_{n+1} \subseteq \C_n$), we get an exact 
sequence of triangulated functors,
\[ \C_{n+1} \rightarrow \C_n \rightarrow \C_n/\C_{n+1} \rar 0.\]
Applying the functor $K_0(-)$ gives an exact sequence of abelian groups,
\[ K_0(\C_{n+1})\rightarrow K_0(\C_n) \rightarrow K_0(\C_n/\C_{n+1}) \rightarrow 0.\]
We do not know much about these groups beyond the fact that they are countably generated abelian groups.
The hard and interesting thing here is to determine these groups and understand the map  
$K_0(\C_{n+1})\rightarrow K_0(\C_n)$ in the above exact sequence. Note that if $X$ is a type-$n$ spectrum and $v$ is a
$v_n$-self map on $X$ \cite{hs}, then $v$ being an even degree self map, the Grothendieck class of cofibre $X/v$ is the zero element in $K_0(\C_n)$. 
Now  heuristically we expect every finite type $n+1$ spectrum to the cofibre of a $v_n$-self map on a type-$n$ spectrum. 
Therefore it is reasonable to conjecture that the map $K_0(\C_{n+1})\rightarrow K_0(\C_n)$ is trivial. 
Also note that this conjecture is equivalent
to the conjecture that $\C_{n+1}$ is contained in every dense triangulated subcategory of $\C_n$. We give some evidence to this
conjecture in subsection \ref{se:lattice}; see proposition \ref{prop:newevidence}.

\subsection{$\C_n$ - Higher thick subcategories.} 
Now we want to study the triangulated subcategories of the thick subcategories $\C_n$, for $n > 1$. 
To this end, we make use of the spectra related to $\BP$ and their homology theories to construct some Euler 
characteristic functions.  Recall that there is a ring spectrum called the Brown-Peterson spectrum (denoted by $\BP$) whose coefficient ring is given by
$\mathbb{Z}_{(p)}[v_1,v_2,\cdots,v_n,\cdots]$, with $|v_i| = 2p^i-2$.
Associated to $\BP$, we have, for each $n \ge 1$, the Johnson-Wilson spectrum
$\bp{n}$ whose coefficient ring is  $\mathbb{Z}_{(p)}[v_1,\cdots,v_n]$.

For each $n\ge1$, recall that there is a generalised Moore spectrum of type-$n$. For $n=1$ this is 
just the Moore spectrum $M(p)$. For $n=2$ this is the cofibre of a self map of the Moore spectrum 
\[ \Sigma^{|v_1^{i_1}|} M(p) \stk{v_1^{i_1}} M(p) \]  
that induces multiplication by $v_1^i$ (for some $i$) in $\BP_*$ homology  and is denoted by
$M(p,v_1^{i_1})$. For higher values this is defined inductively: A type-$n$ generalised Moore spectrum is obtained by taking the
cofibre of a self map
\[\Sigma^{|v_n^{i_n}|} M(p,v_1^{i_1},\cdots,v_{n-1}^{i_{n-1}}) \stk{v_n^{i_n}}
  M(p,v_1^{i_1},\cdots,v_{n-1}^{i_{n-1}}) \]
that induces multiplication by $v_n^{i_n}$ (for some $i_n$) in $\BP_*$ homology and is inductively 
denoted by $ M(p,v_1^{i_1},\cdots,v_{n-1}^{i_{n-1}}, v_n^{i_n})$. For sufficiently large powers of the $v_i$'s these spectra
are known to exist \cite{hs}. However the existential problem of such spectra with specified exponents seems to be a very hard question.

\begin{lemma} Let $X$ be any spectrum in $\C_n$ ($n \ge 1$). Then  $\bp{n-1}_i X$ is always a finite abelian $p$-group, and is
zero for all but finitely many $i$. 
\end{lemma}

\begin{proof} The strategy here is a thick subcategory argument. Say that a finite $p$-local 
spectrum $X$ has the property P if $\bp{n-1}_i X$ satisfies the conditions in the statement of the lemma.
It is straight forward to verify that the property P is generic.
Now by the Hopkins-Smith thick subcategory theorem, we will be done if we can exhibit one
generic type $n$-spectrum for which P holds. To this end, we consider generalised type-$n$ Moore spectrum
$M(p, v_1^{i_1}, \cdots,v_{n-1}^{i_{n-1}})$.
We have, for each $n\ge 1$,
\[ \bp{n-1}_* M(p, v_1^{i_1}, \cdots,v_{n-1}^{i_{n-1}} ) = 
\frac{\mathbb{Z}_{(p)}[v_1,\cdots,v_{n-1}]}{(p,v_1^{i_1},\cdots,v_{n-1}^{i_{n-1}})}.  \]
This homology is a finitely generated $\mathbb{F}_p$ algebra, and therefore the generalised Moore spectrum in question 
has the property P. So we are done.  In fact it is also clear  that 
the the full subcategory of all finite $p$-local spectra that satisfy the property P is precisely
the thick subcategory $\C_n$.
\end{proof}

With the above lemma at hand, we can define a function $\chi_n : \C_n \rightarrow \mathbb{Z}$  as
\begin{equation} \label{E:chi_n}
 \chi_n X = \sum_{i=-\infty}^{\infty} (-1)^i \log_p |\bp{n-1}_i (X) |.
\end{equation}
The previous lemma verifies that this is a well-defined function. It is straight forward to 
verify that this function is an Euler characteristic function. Moreover the $\bp{n-1}_*$ homology of the generalised type-$n$ Moore
spectrum $M(p, v_1^{i_1}, \cdots,v_{n-1}^{i_{n-1}})$
is non-trivial and is concentrated in a finite range of even degrees and hence
$\chi_n M(p,v_1^{i_1}, \cdots, v_{n-1}^{i_{n-1}})$ is non-zero.  So by the universal property
of the Euler characteristic function $\chi_n$, we have the following split short exact sequence.
\[
\xymatrix{
0 \ar[r] & \Ker(\phi_n) \ar[r] & K_0(\C_n) \ar[r] \ar[dr]_{\phi_n} & \im (\phi_n) \ar[r] \ar@{^{(}->}[d] & 0 \\  
         &                     &                  & \ints                             & 
}
\]
This discussion can be summarised in the following proposition.

\begin{prop} For each $n\ge 1$, $K_0(\C_n)$ has a direct summand isomorphic to $\mathbb{Z}$.
 This gives, for each 
$k \ge 0$, a dense triangulated subcategory of $\C_n$ defined by 
\[ \C_n^k := \{ X \in \C_n : \chi_n (X) \equiv 0\!\!\! \mod{l_n k} \} \]
where $l_n$ is a generator for the cyclic group $\im (\phi_n)$.
\end{prop}

Note that these triangulated subcategories correspond to the subgroups $n\mathbb{Z} \oplus \Ker(\phi_n)$ and hence are dense in $\C_n$
by theorem \ref{main}.

\begin{rem} The above proposition  recovers the dense triangulated subcategories $\C_1^k$  of theorem \ref{thm:steve}. In fact, the 
Euler characteristic function (\ref{E:chi_n}) agrees with (\ref{E:chi_1}) when $n=1$. This is because
$\BP \la 0 \ra \cong H\mathbb{Z}_{(p)}$ (the Eilenberg-Mac Lane spectrum for the integers localised at $p$),
and for $p$-torsion spectra, ${H\mathbb{Z}_{(p)}}_* (X) \cong  H\mathbb{Z}_* (X)$. 
Thus, for $X$ in $\C_1$, we have $\bp{0}_* X \cong H\mathbb{Z}_* (X)$ and consequently the two Euler characteristic functions agree.
\end{rem}
 
\subsection{A conjecture of Frank Adams}
It is easily seen that if the Smith-Toda complex $V(n-1) ( := M(p,v_1,v_2,\cdots v_{n-1}))$ exists, then the map $\phi_n: K_0(\C_n) \rar \ints$ is surjective. 
The natural thing to do now is to determine the kernel of $\phi_n$. This turns out to be a very hard problem. It is known that  
$V(1)$ exists at odd primes. In view of this Frank Adams made the following conjecture.

\begin{conj} \label{conj:Adams}
\cite[Page 529]{MayTho} The thick subcategory $\C_2$ (at odd primes) can be generated by iterated cofiberings of the Smith-Toda complex $V(1)$.
\end{conj}

This conjecture is equivalent to saying that $\Ker(\phi_2)=0$, or equivalently that $K_0(\C_2) \cong \mathbb{Z}$. 
Adams \cite[Page 528]{MayTho} also asked the following weaker question.

\vskip 3mm 
\noindent
\textbf{Question:} What is a good set of generators for $\C_n$ ?  (A set $A$ generates $\C_n$ if the smallest 
triangulated category that contains $A$ is $\C_n$.) 
\vskip 3mm

This is a very important question and one answer was given by Kai Xu \cite{xu}. Before we can state his result we
need to set up some terminology.

Recall that a spectrum $X$ is \emph{atomic} if it does not admit any non-trivial idempotents, i.e., if $f \in [X,X]$
is such that $f^2 = f$, then $f = 0$ or $1$. Since idempotents split in the stable homotopy category, this is
also equivalent to saying that $X$ does not have any non-trivial summands. Now if $X$ is a finite $p$-torsion atomic
spectrum, then the finite non-commutative ring $[X, X]$ of degree zero self maps of $X$ 
modulo its Jacobson radical (intersection of all the left maximal ideals)
is isomorphic to a finite field \cite{adamskuhn}:
\[ [X,X]/rad \cong \mathbb{F}_{p^k}, \; \mbox{for some} \; k.\]
So for every $p$-torsion atomic spectrum $X$, we define $e(X)$ to be the integer $k$  given by the above isomorphism. 

\begin{example} It can be easily verified that for all $i \ge 1$, $M(p^i)$ is an atomic spectrum with $e(M(p^i))=1$.
\end{example}

The main result of \cite{xu} which uses the nilpotence results \cite{dhs,hs} as the main tools then states:

\begin{thm}\cite{xu} For any every pair of natural numbers $(n, k)$, there is an atomic spectrum $X$
of type-$n$ such that $e(X)=k$. Further, if $C(n,k)$ denotes the triangulated subcategory of $\C_n$
generated by the type-$n$ atomic spectra with $e(-) \le k$, then
$C(n,k) = \C_n$.
\end{thm}
 
So this gives an answer to Adams's question. One set of generators for $\C_n$ can be taken to be the collection
of all type $n$ atomic spectra $X$ with $e(X)=1$.  The next natural question is how big is this set? Is this a finite set? If so, 
then we can infer that the Grothendieck groups $K_0(\C_n)$ are  all finitely generated abelian groups. But we do not know the answer to this question.

Using this theorem of Xu, we can now revise Adams's conjecture as follows.

\begin{conj} At odd primes, $V(1)$ generates (by iterated cofiberings) all  type-$2$ atomic spectra $X$ such that $e(X)=1$.
\end{conj}

\subsection{The lattice of triangulated subcategories of $\C_0$} \label{se:lattice}

We will now study the lattice of triangulated subcategories of 
finite $p$-local spectra. First recall that the thick subcategories of finite $p$-local
spectra are nested \cite{mit}, i.e., 
\[ \cdots \subsetneq \C_{n+1} \subsetneq \C_n \subsetneq \C_{n-1} \subsetneq \cdots \cdots \subsetneq \C_1 \subsetneq \C_0. \]
Our goal now is to understand how the triangulated subcategories $\C_n^k$ fit in this chromatic chain. We begin with a simple proposition.

\begin{prop} $\C_1$ is contained in every dense triangulated subcategory of $\C_0$, i.e., $\C_1 \subsetneq \C_0^k$ for all $k \ge 0$.
\end{prop}

\begin{proof} If $X$ is in $\C_1$, then its rational homology is trivial and hence
its rational Euler characteristic $\chi_0(X)$ is zero. Therefore
$X$ belongs to every dense triangulated subcategory of $\C_0$. It is easy to see that the spectrum $S \vee \Sigma S$ belongs  to $\C_0^k - \C_1$ 
for all $k \ge 0$, therefore the containment $\C_1 \subseteq \C_0^k$ is proper.
\end{proof}

Motivated by this proposition, we wondered if it is true that
for all non-negative integers $n$ and $k$,
\[ \C_{n+1} \subsetneq \C_n^k  \subseteq \C_n. \]
We now proceed to show that this is indeed the case. For better clarity we separate the cases $n=1$ and $n >1$.

\subsubsection{An Atiyah-Hirzebruch spectral sequence.}

We prove that $\C_2 \subsetneq \C_1^k$ for all $k \ge 0$. 
Our main tool will be an Atiyah-Hirzebruch spectral sequence. We begin with two lemmas -- the first one is an elementary algebraic fact and the second one is a standard topological fact.

\begin{lemma}\label{parallel} If $A:= \cdots \rar  0 \rar A_1 \rar  \cdots \rar A_k \rar 0 \rar \cdots$ is a bounded chain complex of finite $p$-groups,
then $\sum_i (-1)^i \log_p |A_i| $= $\sum_i (-1)^i \log_p |H_i(A)|$.
\end{lemma}

\begin{proof} This is left as an easy exercise to the reader.
\end{proof}

\begin{lemma} \label{le:adams}  The thick subcategory $\C_2$ consists  of all finite $p$-torsion spectra whose complex 
$K$-theory is trivial.
\end{lemma}

\begin{proof} We use a result of Adams \cite{adams} which states that the complex $K$-theory localised at $p$ 
splits as a wedge of suspensions of $E(1)$. 
More precisely, $K_p = \bigvee _{i=0}^{i=p-2} \Sigma^{2i} E(1)$. In particular, $\la K_p \ra = \la E(1)\ra$. With
this at hand, we get the following equalities of Bousfield classes.
\begin{eqnarray*}
\la K \ra &=& \bigvee_p \la K_p \ra \\ 
          &=& \bigvee_p \la E(1)\ra \\
          &=& \bigvee_ p ( \la K(0)\ra \bigvee \la K(1) \ra ) . 
\end{eqnarray*}
The last equality follows from \cite[Theorem 2.1(d)]{rav}.
Now it is clear from these equations that for $X$ finite and $p$ torsion, $K_* X = 0$ if and only if $K(1)_* X = 0$.
\end{proof}

\begin{prop} \label{prop:newevidence}  $\C_2$ is properly contained in every dense triangulated subcategory of $\C_1$, i.e., 
$\C_2 \subsetneq \C_1^k$ for all $k \ge 0$. 
\end{prop}

\begin{proof} Recall that $\C_1^k$ is the collection
of $p$-torsion spectra $X$ for which $\chi_1(X)$ is divisible by $k$. So it is clear that $X \in \C_1^k$ for all $k$ 
if and only if $\chi_1(X)=0$. 
Therefore by lemma \ref{le:adams} we have to show that if $X$ is a finite $p$-torsion spectrum for which
$K_*(X) = 0$, then $ \chi_1(X) := \sum (-1)^i \log_p|H\mathbb{Z}_i X | = 0$.

\begin{figure}
\begin{center} 
\scalebox{.75}{\input{ah.pstex_t}}
\caption{Atiyah-Hirzebruch spectral sequence: $E^3$ - term}
\label{E3}
\end{center}
\end{figure}
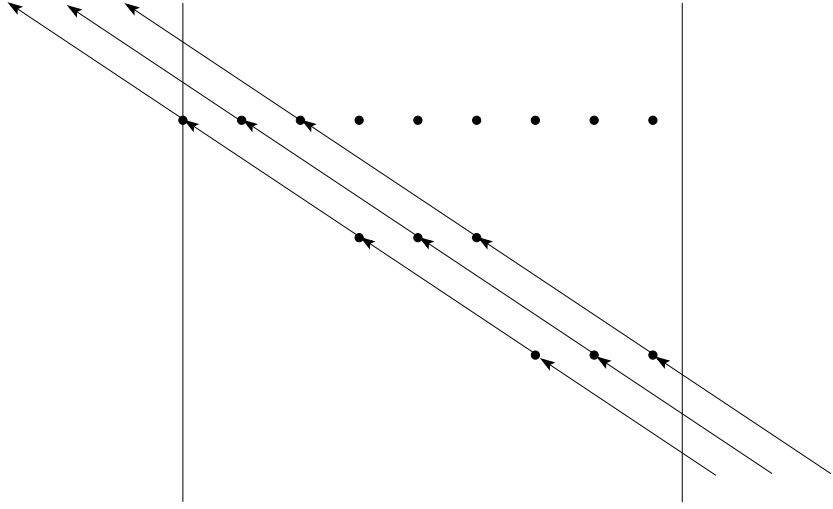

Recall that the spectrum $K$ (also denoted by $\text{BU}$) is a ring spectrum whose 
coefficient group is given by the complex Bott periodicity theorem; $K_* = \mathbb{Z}[u,u^{-1}]$ where $|u|=2$.
We make use of the Atiyah-Hirzebruch spectral sequence 
\[ E^2_{s,t} = H_s(X; K_t) \implies K_{s+t} X. \]
converging strongly to the $K$-theory of $X$. The differentials ($d_r$) in this spectral sequence have bidegrees $|d_r| = (-r, r-1)$.
Note that $X$ being a finite spectrum, the $E^2$ page is concentrated in a vertical strip of finite width (see Figure \ref{E3}) and 
therefore the spectral sequence collapses after a finite stage. Since it converges to $K_*(X)$, which is zero by hypothesis,
we conclude that for all sufficiently large $n$, $E^n = E^{\infty} = 0$.

Next we claim that the function
\[ n \longmapsto \underset{i}{\sum} (-1)^i \log_p |E^n_{i,0}|\]
is a constant function.

Assuming this claim, we will finish the proof of the proposition. When $n=2$, this function 
takes the value  $\sum_i (-1)^i \log_p|H\mathbb{Z}_i\, X | = \chi_1(X)$  and for large enough $n$,
the function takes the value $0$ (since $E^n = E^{\infty} = 0$). Since the function is constant (by the above claim), we get 
$\chi_1(X) = 0$.

Now we prove our claim. First note that this spectral sequence is a module over the coefficient ring $K_* \cong \mathbb{Z}[u,u^{-1}]$. Therefore
the differentials commute with this ring action. Also since $u$ is a unit, it acts isomorphically on the spectral sequence, and therefore
induces periodicity on $E^4$: $E^4_{i,*} \cong E^4_{i-2,*}$. 
Also note that just for degree reasons, all the even differentials
are zero. So at $E_3$, where the first potential non-zero differentials occurs, the alternating sum $\sum_i (-1)^i \text{log}_p |E^3_{i,0}|$
can be broken into three parts as shown in the equation below.
\[ \sum_i (-1)^i \text{log}_p |E^3_{i,0}| = 
   \sum_{i\equiv 0 (3)} +
   \sum_{i\equiv 1 (3)} +
   \sum_{i\equiv 2 (3)}. \]
Now because of the periodicity of the $E^3$ page, we can assemble the terms in  the above equation 
along three parallel lines in the $E_3$ page, where the term with congruence class $l$ modulo three corresponds to the line with $t$-intercept $l$; see figure \ref{E3}.
Now we can apply lemma \ref{parallel} along each of these lines which are bounded 
complexes of finite abelian $p$ groups to pass to the homology groups. Invoking the periodicity of the 
differentials again, we conclude that the new alternating sum thus obtained is
equal to  $\sum_i (-1)^i \text{log}_p |E^4_{i,0}|$. Now an easy induction will complete the proof of the claim:
At ($E^r, d_r$),  for $r$ odd, we decompose the alternating sum $\sum_i (-1)^i \log_p|E_{i,0}^r|$ into $r$ parts, one for each congruence class modulo
$r$, and assemble these terms along $r$ parallel lines on the $E^r$ page such that the term whose congruence class is 
$l$ modulo $r$ corresponds to the line with $t$-intercept $l$. Periodicity of the differentials can be used (as before) 
to complete the induction step.

So we have shown that $\C_2 \subseteq \C_1^k$ for all $k \ge 0$. To see that this inclusion is strict, observe that the 
$p$ torsion spectrum $M(p) \vee \Sigma M(p)$ belongs to $\C_1^k - \C_2$ for all $k$.
\end{proof}

\begin{cor} \label{cor:C1/C2} $K_0(\C_1/C_2) \cong \ints$ generated by the Grothendieck class of the image of the Moore spectrum under the quotient functor
$\C_1 \rar \C_1/\C_2.$
\end{cor}

\begin{proof} On applying $K_0(-)$ to the sequence $\C_2 \rar \C_1 \rar \C_1/\C_2 \rar 0$, we get an exact sequence of abelian groups:
$K_0(\C_2) \rar K_0(\C_1) \rar K_0(\C_1/\C_2) \rar 0$; see (\ref{eq:K_0exactseq}).  The first map in this sequence is the zero map by proposition 
\ref{prop:newevidence}, and $K_0(\C_1)$ was shown to be infinite cyclic on $[M(p)]$. So the corollary follows by combining these two facts.
\end{proof}

The inclusion $\C_2 \subsetneq \C_1^k \; \forall \, k$, which we have just established, gives some new evidence to conjecture \ref{conj:Adams} of Adams. 
To see this, first note that since $M(p,v_1)$ is a cofibre of an even degree self map 
of $M(p)$,  $\chi_1(M(p,v_1)) = 0$. Now if  $M(p,v_1)$ generates $\C_2$  by cofibre sequences (Adam's conjecture), 
then  $\chi_1(X)=0$ for all $X \in \C_2$. This now clearly implies that $\C_2 \subsetneq \C_1^k \; \forall \, k$.

One can try to test this conjecture by asking whether we can generate $M(p^i,v_1^j)$ (whenever it exists) from $M(p,v_1)$.
Results from nilpotence and periodicity \cite{hs} give the following  partial answer to this question.

\begin{prop} For every fixed positive integer $k > 0$, there exists infinitely many positive integers $j$ for which
$M(p^k,v_1^j)$ can be generated from $M(p,v_1)$ using cofibre sequences.
\end{prop}

We leave the proof of this well-known proposition to the reader.

All these results give only some evidence for Adams's conjecture. The conjecture, however, still remains open.

\subsubsection{A Bockstein spectral sequence}

Our goal now is to prove: $\C_{n+1} \subsetneq \C_n^k$ for all $n$ and $k$. We mimic our strategy for the case
$n=1$ by replacing the Atiyah-Hirzebruch spectral sequence with a Bockstein spectral sequence.  The above inclusion 
is an easy corollary of the following theorem.

\begin{thm} If $X$ is spectrum in $\C_{n+1}$, then  
\[\chi_n (X) := \sum (-1)^i \log_p |\bp{n-1}_i X| = 0.\]
\end{thm}

\begin{cor} $\C_{n+1} \subsetneq \C_n^k$ for all $k \ge 0$ and all $n \ge 1$.
\end{cor}

\begin{proof} By the above theorem, for a spectrum $X$ in $\C_{n+1}$, $\chi_n(X) = 0$. Therefore $X$ belongs to $\C_n^k$ for all $k$.
The spectrum $F \vee \Sigma F$, where $F$ is a  type-$n$ spectrum, belongs to $\C_n^k - \C_{n+1}$. So $\C_{n+1} \subsetneq \C_n^k$ for all $k \ge 0$.
\end{proof}

\n
We now outline the strategy for proving the above theorem. This is very similar to the proof of proposition \ref{prop:newevidence}.
First recall that $\C_{n+1}$ can also be characterised as 
the collection of finite $p$-local spectra that are acyclic with respect to $E(n)$. So we seek a strongly convergent spectral sequence 
                                 \[ E_r^{**} \implies E(n)_* X,\] 
whose $E_1$ term is build out of  $\bp{n-1}_* X$. We show that a certain Bockstein spectral sequence has this property
and that it collapses after a finite stage. Finally we work backward to conclude that
$\chi_n X = 0$. 

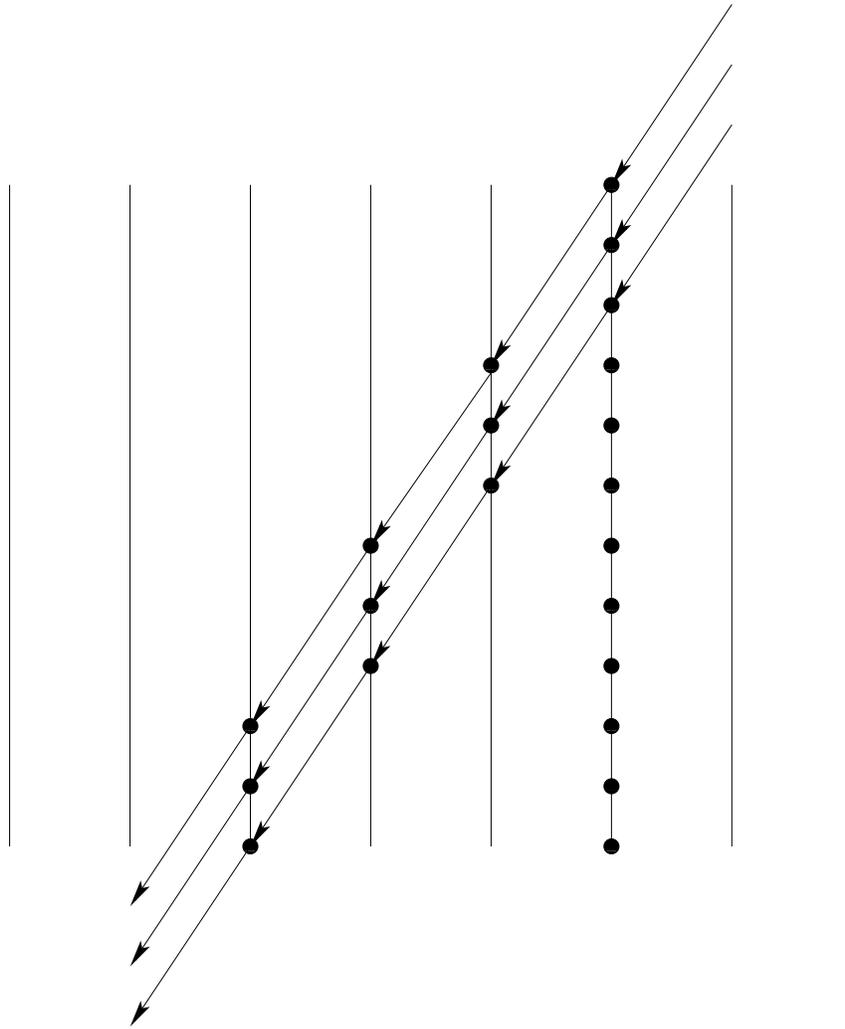
\begin{figure}
\begin{center}
\scalebox{.75}{\input{bss.pstex_t}}
\caption{Bockstein spectral sequence: $E_1$ - term}
\label{E1}
\end{center}
\end{figure}

Now we proceed to construct such a spectral sequence. Fix an integer $n \ge 1$ and recall that 
\[E(n) = v_n^{-1} \bp{n} = \hocolim{v_n} \bp{n}.\]
This gives a sequence of cofibre sequences that fit into a diagram extending to infinity in both directions as shown below.
\[
\xymatrix@=1.1em{
\cdots \ar[r]^{v_n\ \ \ \ \ \ } & \sus{|v_n|} \bp{n} \ar[r]^{\ \ v_n} \ar[d] & \bp{n} \ar[r]^{v_n \ \ \ \  } \ar[d] & \sus{-|v_n|} \bp{n} \ar[r]^{\ \ \ \ v_n} \ar[d] 
& \cdots \cdots \ar[r] & E(n) \\
\cdots              &  \sus{|v_n|} \bp{n-1} \ar[ul]|{\circ}           &        \bp{n-1} \ar[ul]|{\circ}  &  \sus{-|v_n|} \bp{n-1} \ar[ul]|{\circ}  
& \cdots  &                 
}
\]
Since both the functors  $(-)\wedge X$ and $\pi_*(-)$ commute with the functor $\hocolim{v_n}(-)$, smashing the above
diagram with $X$ and taking $\pi_*$ gives an exact couple of graded abelian groups:
\[
\xymatrix@=0.65em{
\cdots \ar[r]^{v_n \ \ \ \ \ \ \ \ } & \sus{|v_n|} \bp{n}_*X \ar[r]^{\ \ \ \ v_n} \ar[d] & \bp{n}_*X \ar[r]^{v_n\ \ \ \ \ \ \  } \ar[d] & \sus{-|v_n|} \bp{n}_*X  \ar[d] 
 \cdots \cdots \ar[r] & E(n)_*X \\
\cdots              &  \sus{|v_n|} \bp{n-1}_*X \ar[ul]|{\circ}          &        \bp{n-1}_*X \ar[ul]|{\circ}  &  \sus{-|v_n|} 
\bp{n-1}_*X \ar[ul]|{\circ}  &                    
}
\]
This exact couple gives rise to a (Bockstein) spectral sequence $E_r^{*,*}$ in the usual way. We choose a convenient grading so that 
the $E_1$ term is concentrated in a horizontal strip of finite width, i.e., $E_1^{*,q} = 0$ for $|q| >> 0$ (see figure \ref{E1}). This can be ensured by setting
\[D_1^{p,q} = \Sigma^{-p|v_n|} \bp{n}_{-p|v_n|+q} X, \]
\[E_1^{p,q} = \Sigma^{-p|v_n|} \bp{n-1}_{-p|v_n|+q} X.\]
With this grading one can easily verify that the differentials have bidegrees given by $|d_r|=(-r, -r|v_n|-1)$. 
It is clear from figure \ref{E1} that after a finite stage all the differentials
exit the horizontal strip and therefore the spectral sequence collapses. Another important fact that we need about this
Bockstein spectral sequence is the periodicity of all the differentials. More precisely, for all integers $p$ and $r \ge 0$,
$E_r^{p,*} = E_r^{p+1,*}$, and further the following diagram commutes.
\[
\xymatrix{
E_r^{p,*}\;\;\;  \ar[d]^{=}  \ar[r]^{d_r\;\;\;\;\;\;\;} & \;\;E_r^{p-r,*-1-r|v_n|} \ar[d]^{=}\\
E_r^{p+1,*}  \ar[r]^{d_r\;\;\;\;\;\;\;\;} & \;  E_r^{p+1-r,*-1-r|v_n|}
}
\]
Now we show that this spectral sequence converges strongly to $E(n)_*(X)$. For this, we make use of a theorem
of Boardman \cite{boardman}. Before we can state his result we have to recall some of his terminology. Consider an exact couple 
of graded abelian groups:

\[
\xymatrix{
\cdots \ar[r]^{i} & A^{s+1} \ar[r]^{i} \ar[d]^j & A^s \ar[r]^{i}\ar[d]^j & A^{s-1} \ar[r]^{i} \ar[d]^j 
& A^{s-2} \ar[r]^{i} \ar[d]^j & \cdots  \ar[r] & A^{-\infty} \\
\cdots              &  E^{s+1} \ar[ul]^k|{\circ}           &   E^s \ar[ul]^k|{\circ}  &  E^{s-1} \ar[ul]^k|{\circ}  
& E^{s-2} \ar[ul]^k|{\circ}  &\cdots \cdots  &                 
}
\]

This gives the following filtration of the groups $E^s$ by cycles and boundaries of the differentials
in the spectral sequence that arises from this exact couple:
\[ 0 = B_1^s \subset B_2^s \subset B_3^s \subset \cdots \
   \cdots \subset Z_3^s \subset Z_2^s \subset Z_1^s=E^s, \]
where $Z_r^s := k^{-1}(\text{Im}\, [i^{(r-1)}: A^{s+r} \rar A^{s+1}])$, and $B_r^s:= j \ker[i^{(r-1)}: A^s \rar A^{s-r+1}]$.

Associated to this exact couple and the resulting spectral sequence, \cite{boardman} defines the following groups.
These definitions also hold when $n=\infty$.
\begin{itemize}
\item $A^{-\infty}:= \colim{s} A^s$, \hspace{3 mm}  $A^{\infty}:= \clim{s} A^s$, \hspace{3 mm} $RA^{\infty}:= \rlim{s} A^s$
\item $K_n A^s:= \ker[i^{(n)}: A^s \rar A^{s-n}] $, \hspace{3mm} $\text{Im}^r A^s:= \text{Im}[i^{(r)}: A^{s+r} \rar A^s ]$
\item $K_n \text{Im}^r A^s := K_n A^s \cap \text{Im}^r A^s$
\item $W=\colim{s} \rlim{r} K_{\infty} \text{Im}^r A^s $
\item $RE_{\infty}^s = \rlim{r} Z_r^s$.
\end{itemize}

The main theorem of \cite{boardman} then states:

\begin{thm}\cite[Theorem 8.10]{boardman}
The spectral sequence arising from the above exact couple converges strongly to $A^{-\infty}$ if the          
obstruction groups $A^{\infty}$, $RA^{\infty}$, $W$, and $RE_{\infty}$ are zero.
\end{thm}

We recall a few elementary facts about inverse limits before we can apply this theorem to our Bockstein 
spectral sequence. The proofs follow quite easily from the universal properties of these
limit functors. Parts of this lemma can also be derived from the more general Mittag-Leffler condition.

\begin{lemma} \label{facts-limits} Consider a sequence of groups (graded) and homomorphisms:
\[ \cdots \rar A^{s+1} \rar A^s \rar A^{s-1} \rar \cdots.\]
Then the following statements hold.
\begin{itemize}
\item If, for some integer $M$, the composite of $M$ consecutive maps in this sequence is always null, 
then  $\clim{s} A^s = \rlim{s} A^s = 0$.
\item If there is an integer $L$ such that for all $s \ge L$ the map $A^{s+1} \rar A^s$ is the identity
map, then $\clim{s} A^s = A^{L}$ and $\rlim{s} A^s = 0$
\end{itemize}
\end{lemma}

We now show that our Bockstein spectral sequence converges strongly to $E(n)_* X$ by verifying the hypothesis of
Boardman's theorem. So we apply his theorem to the sequence,
\begin{equation}\label{vnseq}
\cdots \stk{v_n}  \sus{|v_n|} \bp{n} \stk{v_n} \bp{n} \stk{v_n} \sus{-|v_n|} \bp{n} \stk{v_n} \cdots  
\end{equation}
For brevity, we shall denote $\sus{-s|v_n|} \bp{n}_* X$ by $A^s$.

\n
(a) $A^{\infty}=0$: By lemma \ref{facts-limits}, all we need to show is that there is some integer $M$
such that the composite of any $M$ consecutive maps in the sequence \eqref{vnseq}
is zero. $X$ being a spectrum in $\C_{n+1}$, $\bp{n}_*X$ is concentrated only
in a finite range. Now since $v_n$ is a graded map of degree $2(p^n-1)$, a sufficiently large iterate of $v_n$
vanishes, hence $A^{\infty} = \clim{s} \sus{-p|v_n|} \bp{n}_* X = 0$.

\n
(b) $RA^{\infty}=0$: This is also immediate from part (1) of lemma \ref{facts-limits} because we have already seen that a sufficiently large
iterate of $v_n$ is zero in part (a).

\n
(c) $W=0$:  The colimit of the sequence \eqref{vnseq} is $E(n)_*X$ and this latter group is zero by hypothesis.
It now follows that $K_{\infty} A^s = A^s$. We make use of this fact in the third equality below. 
\begin{eqnarray*}
W &=& \colim{s} \rlim{r} \left( K_{\infty} \text{Im}^r A^s \right) \\
  &=& \colim{s} \rlim{r} \left( K_{\infty} A^s \cap \text{Im}^r A^s \right) \\
  &=& \colim{s} \rlim{r} \left( A^s \cap \text{Im}^r A^s \right) \\
  &=& \colim{s} \rlim{r} \; \text{Im}^r A^s \\
  &=& \colim{s} \rlim{r} ( \cdots 0 \rar 0 \rar \cdots \subseteq \text{Im}^2 A^s \subseteq \text{Im}^1 A^s) \\
  &=& \colim{s} 0  = 0 
\end{eqnarray*}
\n
(d) $RE_{\infty}=0$: This follows from the fact that the spectral sequence collapses after a finite stage. For each fixed $s$, 
all inclusions in the sequence 
\[ \cdots \stk{=} Z_{r+1}^s \stk{=} Z_{r}^s \subseteq \cdots \subseteq Z_3^s \subseteq Z_2^s \subseteq Z_1^s = E_1^s ,  \]
become equalities eventually. So we now invoke part (2) of  lemma \ref{facts-limits} to conclude that 
$RE_{\infty}^s = \rlim{r} Z_r^s = 0$. \\

So we have shown that all the obstruction groups vanish and therefore Boardman's theorem tells us that our
spectral sequence converges strongly to $E(n)_*X$. \\

We now move to the final part of the theorem: $\chi_n(X)=0$. We mimic the argument given for
the fact ``$\C_2 \subseteq \C_1^k$''. The crux of the proof lies in the key observation that the $q$ component
of the differential $d_r$ is always an odd number ($-r|v_n|-1$). We claim that the function
\[k \overset{\phi}{\longmapsto} \sum(-1)^i \log_p |E_k^{0,i}|\]
is constant. Toward this, we decompose the
alternating sum  $\sum(-1)^i \log_p |E_1^{0,i}|$ into $t:=|v_n|+1$ parts as follows:
\[ \sum_{i} (-1)^i \log_p |E_1^{i,0}| = \sum_{i\equiv 0 (t)} + \sum_{i \equiv 1 (t)} + \cdots + \sum_{i \equiv t-1 (t)}\] 
Now the periodicity of all the differentials coupled with the fact that $|v_n|+1$ is an odd number will enable us to  assemble all these terms 
on the right hand side of this equation
along $t$ parallel complexes of differentials on the $E_1$ term. (The term corresponding to the
congruence class $l$ modulo $t$ will correspond to the parallel complex with $q$ intercept $l$; see figure \ref{E1}.) 
We can now apply lemma \ref{parallel} to each of these complexes and pass on to the homology groups
without changing the underlying alternating sum. 
Again using the periodicity of the differentials, we can reassemble all the terms after taking homology to obtain $\phi(2)$. This shows that
$\phi(1)=\phi(2)$. Now a straight forward induction will complete the proof of the claim.

Finally, to see that $\chi_n(X)=0$, observe that when $k=1$, $\phi$ takes the value $\chi_n(X)$, and for all sufficiently
large values of $k$, $\phi$ is zero because our spectral sequence collapses at a finite stage and  converges 
strongly to $E(n)_*X$, which is zero by hypothesis. Since we know that $\phi$ is constant, this completes the
proof of the theorem.

\subsection{Questions}

\subsubsection{Grothendieck groups and Adams's Conjecture}
The biggest problem that needs to settled is the classification of all triangulated subcategories of $\F_p$. We believe that such a 
classification  will reveal some hidden ideas behind the chromatic tower which might shed some new light on the stable homotopy category.  
We have seen that $K_0(\C_0)$ is infinite cyclic with the sphere as the generator, and $K_0(\C_1)$ is infinite cyclic with the mod-$p$ 
Moore spectrum as the generator. Computing the Grothendieck groups of higher thick subcategories ($\C_n$, $n \ge 2$), as we have seen, 
is much more complicated. In particular Adams' conjecture remains open.

Another interesting and related question at this point is the following. In $\C_n$, what is the triangulated
subcategory generated by atomic spectra $X$ of type-$n$ for which $e(X)=k$ ($k$ some fixed positive integer)? By \cite{xu}, this is the 
whole of $\C_n$ if $k=1$; otherwise, this will be some dense triangulated subcategory of $\C_n$.  The immediate question that springs out 
now is whether these triangulated subcategories, when $n=2$, are precisely the subcategories $\C_2^k$? 
A non-affirmative answer to this question will settle Adams's conjecture in the negative. Similarly when $n = 1$, it will be interesting 
to match these dense triangulated subcategories with the subcategories $\C_1^k$.

One can also try to investigate some properties of these Grothendieck groups. For instance, are they
finitely generated? are they torsion-free?  

\subsubsection{Euler characteristics and the lattice of triangulated subcategories}

It would be interesting to find an Euler characteristic defined on $\C_n$ that is not a multiple of $\chi_n$ 
(see equation (\ref{E:chi_n})). Such an Euler characteristic function 
might tell something new about $K_0(\C_n)$.

Recall that $l_n$ was defined to the generator of the image of $\phi_n: K_0(\C_n) \rar \ints$. It is clear that $l_n = 1$ if $V(n)$ exists. In general, $l_n$
can be a very large integer and not much is known. 
The following is a conjecture of Ravenel (personal communication): If $(p, f) \ne (2, 1)$, then $p^f$ divides 
$l_{(p-1)f+1}$.

We have seen that the triangulated subcategories $\C_n^k$ are sandwiched between $\C_{n+1}$ and $\C_n$. Is the same true for all triangulated subcategories 
of $\F_p$? 

These are some questions which we think merit further study in this direction.

\section{Triangulated subcategories of perfect complexes} \label{se:classifications} 

For the rest of this paper we work in the derived categories of rings. The subsections that follow are organised as follows.
We begin with a quick recap of the derived category in the next subsection. In subsection \ref{se:algktheory} we review some basic algebraic $K$-theory and
connect it to the problem of classifying triangulated subcategories. We then start applying Thomason's theorem \ref{main} to classify triangulated subcategories of
perfect complexes over PIDs, Artin rings, and some non-noetherian rings. We end with some questions.

Unless stated otherwise, all rings will be assumed to be commutative with a unit.

\subsection{The derived category}

There are many beautiful constructions of the derived category; see \cite{wei} for the classical approach, or \cite{hoveymodel} for a
model category theoretic approach.  We briefly review some preliminaries on the derived category $D(R)$ of a commutative ring $R$. 
It is obtained from the category of unbounded chain complexes of $R$-modules and chain maps 
by  inverting the \emph{quasi-isomorphisms} (maps that induce an isomorphism in 
homology). $D(R)$ is a tensor triangulated category with the derived
tensor product as the smash product and the ring $R$ (in degree $0$) as the unit object. 
It is a standard fact that the small objects of  $D(R)$ are precisely those complexes that are quasi-isomorphic to 
\emph{perfect complexes} (bounded chain complexes of finitely 
generated projective $R$-modules); see \cite[Prop. 9.6]{ch} for a nice proof of this fact. 
It follows from \cite[Corollary 10.4.7]{wei} that the full subcategory of small objects in $D(R)$ is 
equivalent (as a triangulated category) to the chain homotopy category of perfect complexes.
The latter will be denoted by $D^b(\proj\,R)$ and it provides a nice framework for studying small objects.
The full subcategory of small objects in $D(R)$ can also be characterised as the thick subcategory generated by $R$; see \cite[Prop. 9.6]{ch}.

\subsection{Algebraic $K$-theory of rings} \label{se:algktheory}

Now we recall some classical algebraic $K$-theory of rings and connect it to the problem of classifying (dense) triangulated subcategories
of perfect complexes. 
If $R$ is any commutative ring, the \emph{$K$-group of the ring $R$}, which is denoted by $K_0(R)$, is defined to be the free abelian
group on isomorphism classes of finitely generated projective modules modulo the  
subgroup generated by the relations $[P]-[P \oplus Q] - [Q] = 0$, where $P$ and $Q$
are finitely generated projective $R$-modules.
A folklore result says that these $K$-groups are isomorphic to the Grothendieck groups of $D^b(\proj\,R)$.

\begin{prop} (well-known) \label{prop:well-known} If $R$ is any ring (not necessarily commutative), then there is a natural isomorphism of abelian groups
   \[ K_0(R)\cong K_0(D^b(\proj\,R)). \]
\end{prop}
\noindent
\begin{proof}(sketch) Let $A$ denote the free abelian group on the isomorphism classes of perfect complexes 
in the derived category of $R$ and let $B$ denote the free abelian group on isomorphism classes of
finitely generated projective $R$-modules. Since the Grothendieck groups under consideration are quotients of these
free groups, we define maps on $A$ and $B$ that descend to give the desired bijections.
Define $f: A \rightarrow B$ by $f(\la X \ra) = \sum_{i \in \mathbb{Z}} (-1)^i \la X_i \ra$ and
 $g: B \rightarrow A$ by $g(\la M \ra) = \la M[0]\ra$. Now one can easily verify that
these maps descend to the Grothendieck groups and that the descended maps are inverses of each other.
\end{proof}

The importance of this folklore result, for our purpose, can be seen from the following observation.

\begin{rem} This folklore result, along with the theorem \ref{main} of Thomason, connects the 
the problem of classifying dense triangulated subcategories in $D^b(\proj\,R)$ with the
$K$-theory of $R$. More precisely, there is a 1-1 correspondence between the subgroups of $K_0(R)$ and 
the dense triangulated subcategories of $D^b(\proj\,R)$. So this leads us naturally to algebraic $K$-theory -- a subject that has
been extensively studied.
\end{rem}

For the remainder of this section, we review some well-known computations of $K$-groups of rings that will be relevant to us.

\begin{example} \label{ex:K-groups}
If $R$ is any commutative local ring, then it is a fact that every finitely generated projective $R$-module is free. Thus
the monoid $\proj\,R$ (the category of finitely generated projective $R$-modules) is equivalent to the  monoid consisting of the whole numbers.
The Grothendieck group of the later is clearly isomorphic to $\mathbb{Z}$. Thus for local rings, 
\[ K_0(D^b(\proj \, R) \cong K_0(R) \cong \mathbb{Z}. \]
Similarly, it is easy it to see that if $R$ is any principal ideal domain, then $K_0(R) \cong \ints.$
\end{example}

Now we state some results on the $K$-groups of some low  dimensional commutative rings.
These results are of interest because they connect $K$-groups and classical Picard groups of rings.
Before we state the result, we need to recall the definition of the Picard groups of rings.

\begin{defn} The \emph{Picard group}, $\Pic(R)$ of a ring $R$, is defined to be the group of isomorphism classes 
of $R$-modules 
that are invertible under the tensor product. Similarly the \emph{Picard group} of $D(R)$, denoted by $\Pic(D(R))$, 
is the group of isomorphism classes of objects in $D(R)$ that are invertible under the derived tensor product. In both 
these cases, note that the ring $R$ acts as the identity element.
\end{defn}

\begin{thm}\cite{wei1} \label{thm:K_0picard} Let $[\Spec(R), \mathbb{Z}]$ denote the additive group of continuous functions from $\Spec(R)$ to the ring of integers
with the discrete topology. Then the following holds:

\begin{enumerate}
\item For every $0$-dimensional ring $R$, $K_0(D^b(\proj\,R)) \cong [\Spec(R), \mathbb{Z}].$

\item For every $1$-dimensional noetherian ring, 
\[K_0(D^b(\proj\,R)) \cong \Pic(D(R)).\]
\end{enumerate}
\end{thm}

\begin{proof} The first statement is Pierce's theorem; see \cite[Theorem 2.2.2]{wei1}. The second statement can be seen as a corollary of a theorem due to  
Fausk \cite{fau}: There is a natural split short exact sequence (for any commutative ring),
\[ 0 \rightarrow \Pic(R) \rightarrow \Pic(D(R)) \rightarrow [\Spec(R), \mathbb{Z}] \rightarrow 0.\]
Therefore $\Pic(D(R)) \cong \Pic(R) \oplus [\Spec(R),\mathbb{Z}].$ It is shown in \cite[Corollary 2.6.3]{wei1}
that 
\[K_0(D^b(\proj\,R) \cong  \Pic(R) \oplus [\Spec(R),\mathbb{Z}].\] 
So the second statement follows by combining these two results.
\end{proof}

With these and related results from algebraic $K$-theory, one can compute the $K$-groups of various families of rings and that will help
understand the dense triangulated subcategories of perfect complexes over all those rings. However, in order to classify all 
triangulated subcategories of perfect complexes, one has to compute the Grothendieck groups of the thick subcategories of these 
complexes.

We now recall some definitions and a theorem due to Hopkins and Neeman \cite{Ne}
which classifies the thick subcategories of perfect complexes over a noetherian ring.

\begin{defn} \cite{Ne} Given a perfect complex $X$ in $D(R)$, define the \emph{support} of $X$, denoted by 
$\Supp(X)$,  to be the set $\{ p \in \Spec(R) : X \otimes R_p \ne 0 \}$, where $R_p$ is the localisation of 
$R$.
\end{defn}

\begin{defn}
A subset of $\Spec(R)$ is said to be \emph{closed under specialisation} if it
is a union of closed sets under the Zariski topology. Equivalently, and more explicitly, 
a subset $S$ of $\Spec(R)$ is specialisation closed if whenever a prime ideal $p$ is in $S$, then
so is every prime ideal $q$ that contains $p$.
\end{defn}

Now we are ready to state the celebrated thick subcategory theorem of Hopkins and Neeman.

\begin{thm}\cite{Ne} If $R$ is any noetherian ring, then there is a natural order preserving bijection between
the sets,
\begin{center}
\{thick subcategories $\A$ of $D^b(\proj \, R)$\}
\begin{center}
$f\downarrow \ \ \uparrow g$
\end{center}
\{subsets $S$ of $\Spec(R)$ that are  closed under specialisation\}.
\end{center}
The map $f$ sends a thick subcategory $\A$ to  $\bigcup_{X \in \A} \Supp(X)$,
and the map $g$ sends a specialisation-closed subset $S$ to the thick subcategory  $ \T_S:= \{ X \in D^b(\proj \, R): \Supp(X) \in S \}$. 
\end{thm}

The following corollary is immediate from the above theorem.
\begin{cor}  Every thick subcategory of $D^b(\proj\,R)$ is a thick ideal.
\end{cor}

\begin{cor} Let $R$ be any commutative ring such that $K_0(R) \cong \ints$. Then every triangulated subcategory of $D^b(\proj\,R)$ is a triangulated ideal.
\end{cor}

\begin{proof} Since we know that the thick subcategories of $D^b(\proj\,R)$ are thick ideals, they can be viewed as triangulated modules over $D^b(\proj \,R)$. Now we can
apply theorem \ref{refree-main} to conclude that every triangulated subcategory is also a triangulated submodule because $K_0(D^b(\proj\,R)) \cong \ints.$
\end{proof}

It is clear that the proofs of the above corollaries generalise to prove the following proposition.
\begin{prop} Let $\T$ be a tensor triangulated category. If the unit object $S$ is small and generates $\T$, then every thick subcategory 
of compact objects is a thick ideal. Moreover, if the Grothendieck ring of the compact objects in $\T$ is isomorphic to $\ints$, then every triangulated 
subcategory of compact objects is a triangulated ideal. 
\end{prop}

We will now apply Thomason's recipe to classify the triangulated subcategories of perfect complexes over some commutative
rings: Principal ideal domains, Artin rings, and non-noetherian rings with a unique prime ideal.

\noindent
\subsection{Principal ideal domains}

We first set up a few notations and recall some definitions and basic facts about PIDs.

For any element $x$ in $R$, we define the mod-$x$ Moore complex $M(x)$, in analogy with the stable homotopy
category of spectra, to be the cofibre of the self map (of degree 0)
\[  R \stk{x} R \] 
in $D(R)$.

Now we recall the notion of the length of a module. For an $R$-module $M$, a chain
\[M = M_0 \supsetneq M_1 \supsetneq M_2 \supsetneq \cdots \supsetneq M_r = 0\]
is called a composition series if each $M_i/M_{i+1}$ is a simple module (one that does not have any non-trivial submodules).
The length of a module, denoted by $l(M)$, is defined to be the length of any composition series
of $M$. The fact that this is well-defined is part of the Jordan-Holder theorem.

The function $l(-)$ is an additive function on the subcategory of $R$-modules which have finite length, i.e., if 
\[ 0 \rightarrow M_1 \rightarrow \cdots \rightarrow M_s \rightarrow 0\]
is an exact sequence of $R$-modules of finite length, then
\[ \sum_i (-1)^i l(M_i) = 0.\]
Also note that when $R$ is a PID, every finitely generated torsion module has
finite length. (This can be seen as an immediate consequence of the structure theorem for finitely generated 
modules over a PID.)

Finally, we need the following easy exercise. If $p$ and $q$ are two distinct (nonzero) prime elements in a PID $R$, then
for any $i \ge 1$,
\[  
R/(p^i) \otimes_R R_{(q)} =
 \left\{
   \begin{array}{ll}
      R/(p^i) & \mbox{when $p=q$},\\
      0       & \mbox{when $p \neq q$}.
   \end{array}
 \right.
\]

Now we are ready to compute the Grothendieck groups for thick subcategories  of perfect complexes over a PID.
Given a subset $S$ of $\Spec(R)$ that is closed under specialisation, the thick subcategory that 
corresponds to the subset $S$ (under the Hopkins-Neeman bijection) will be denoted by $\T_S$.

\begin{thm} \label{thm:K_0PID}
Let $R$ be a PID and let $S$ be a specialisation closed subset of $\Spec(R)$. Then we have the following.
\begin{enumerate}
\item If $S = \Spec(R)$, then $K_0(\T_S)$ is an infinite cyclic group generated by $R$. 
\item  If $S \ne \Spec(R)$, then $K_0(\T_S)$ is a free abelian group on  the Moore complexes $M(p)$, for 
$p \in S$.
\end{enumerate}
\end{thm}

\begin{proof} 
The first part follows from the fact that every finitely generated projective module over a PID is free;
consequently its Grothendieck group is an infinite cyclic group (see proposition \ref{prop:well-known} and example \ref{ex:K-groups}).
For the second part, first note that a specialisation closed subset $S \ne \Spec(R)$ is a subset of maximal ideals in $R$ 
(because non-zero prime ideals in a PID are also maximal). For each prime element $p$ in $S$, define an Euler characteristic function
$\lambda_p : \T_S  \rightarrow \mathbb{Z}$ by 
\[\lambda_p(X) := \sum_i (-1)^i \; l[H_i(X) \otimes_R R_{(p)}].\]
(Since $S$ does not contain $(0)$, it follows that $H_*(X)$ is a torsion $R$-module. Therefore $\lambda_p(-)$ is a well-defined
Euler characteristic function.)

Also, since $\lambda_p(M(q)) = \delta_p^q$,  the Euler characteristic map 
\[\bigoplus_{p \in S} \lambda_p :  K_0(\T_S) \rightarrow \bigoplus_{p \in S} \mathbb{Z}\]
is clearly surjective.

To see that this map is injective, it suffices to show that every complex
in $\T_S$ can be generated by the set $\{M(p) : p \in S \}$ using cofibre sequences. We do this by induction 
on $\sum_i l[H_i(-)]$. 
If $X \in \T_S$ is such that $\sum_i l[H_i(X)] = 1$, then there exists an integer $j$ such that 
$H_i(X) = 0$ for all $i \neq j$, and $H_j(X) = R/(p)$ for some prime  $p$ in $R$ (because every 
simple module over a PID is of the form $R/(p)$ for some prime  $p$). 
Such an $X$ is clearly quasi-isomorphic to $M(p)$. Now consider an $X \in \T_S$ for which  $\sum_i l[H_i(X] > 1$.
Without loss of generality we can assume that up to suspension $X$ is of the form 
\[ \cdots 0 \rightarrow P_k \rightarrow \cdots \rar  P_1 \rightarrow P_0 \rightarrow 0 \cdots\]
with $H_0(X) \neq 0$ (otherwise, we can replace $X$ with a quasi-isomorphic complex in $\T_S$ which has 
this property). Now pick a non-zero element in $H_0 (X)$
and represent it with a cycle $t$. Since $H_0(X)$ is a torsion module, there exists a prime $p$ and a positive integer $k$ such that $p^k t=0 $ in homology. Replacing 
$t$ with $p^{k-1}t$, we can assume that $pt=0$ in homology, which  means $pt$ is a boundary. So there is an element $y \in P_1$ which maps under the differential
to $pt$. Consider figure \ref{fig:killingahomologyclass}
\begin{figure}
\[
\xymatrix{ \vdots \ar[d] \ar[r] & \vdots \ar[d] \\
           0 \ar[d] \ar[r] & P_2  \ar[d]       \\
           R \ar[d]_p \ar[r]^a & P_1  \ar[d]^c \\
           R \ar[d] \ar[r]^b & P_0  \ar[d]     \\
           0  \ar[r] & 0        
           }
\]
\caption{Killing a homology class in the Hurewicz dimension}
\label{fig:killingahomologyclass}
\end{figure}
where $a(1)= y$, $c(y) = pt$, and $b(1)=t$. This diagram shows a chain map 
between the two complexes in $\T_S$ such that the induced map in homology in dimension $0$ 
sends $1$ to $t$ (by construction). So the class $t$ is killed. Now if we extend this morphism to a triangle
 ($M(p) \rightarrow X \stk{d} Y \rightarrow \Sigma M(p)$)
and look at the long exact sequence in homology, it is clear that 
$X$ and $Y$ have the same homology in all dimensions except dimension 0.
In dimension $0$, part of the long exact sequence gives a short exact 
sequence $ 0 \rightarrow R/(p) \rightarrow H_0(X) \rightarrow H_0(Y) \rightarrow 0$.
Therefore $l[H_0(Y)] = l[H_0(X)] -1$. By induction hypothesis we know that 
$Y$ can be generated by the set $\{M(p): p \in S \}$ using cofibre sequences. The above exact triangle then
tells us that $X$ can also be generated using cofibre sequences in this way.
So we are done.
\end{proof}

\begin{cor} If $R$ is any PID, then $\Pic(D(R)) \cong \ints.$
\end{cor}
\begin{proof} Since $K_0(D^b(\proj\,R)) \cong \ints$, the corollary follows by invoking  theorem \ref{thm:K_0picard}.
\end{proof}

\noindent
\emph{Classification of the triangulated subcategories of $D^b(\proj\,R)$ when $R$ is a PID.} There are two families of triangulated subcategories:
triangulated subcategories that correspond to $S = \Spec(R)$ and  the ones that correspond to $S \ne \Spec(R)$ (subsets of maximal ideals). \\
\noindent 
1.  $S = \Spec(R)$:   Consider the Euler characteristic function 
\[\chi(X) = \sum_{-\infty}^{\infty} (-1)^i \dim_F \{H_i(X) \otimes_R F\},\]
where $F$ is the field of fractions of our domain $R$. For every integer $k$, we define 
\[ D_k = \{ X : \chi(X) \equiv 0 \! \!\mod{k} \}. \]
These are all the triangulated subcategories that are dense in $D^b(\proj\,R)$.\\
\noindent 
2. $S \ne \Spec(R)$:  Given such a subset $S$
   and a subgroup $H$ of $\bigoplus_{p \in S} \mathbb{Z}$, we define 
\[\T(S,H) = \{X \in \T_S: \left(\oplus_{p \in S} \lambda_p \right)(X) \in H \}. \] 
These are all the triangulated subcategories that are dense in $\T_S$.

It is clear from theorem \ref{thm:K_0PID} and theorem \ref{main} that every triangulated subcategory of $D^b(\proj\,R)$ is one of these
two types. 

Here is an interesting consequence of the above theorem. 

\begin{cor} \label{cor:questionBPID} Let $X$ and $Y$ be perfect complexes over a PID. Then $Y$ can be generated from
$X$ using cofibrations if and only if
\begin{itemize}
\item $\Supp(Y) \subseteq \Supp(X)$,  and
\item  If $(0) \in \Supp(X)$, then $\lambda_0(X)$ divides $\lambda_0(Y)$; otherwise, 
  $\lambda_p(X)$ divides $\lambda_p(Y)$ for all $p \in \Supp(X)$.
\end{itemize}
\end{cor}

\subsection{Product of rings: Artin rings} We now  address the following question.\\
\noindent
\textbf{Question:} Suppose a commutative ring $R$ is a direct product of rings. Let us say
\[R \cong R_1 \times R_2 \times \cdots \times R_k,\]
and suppose that we have a classification of all the triangulated subcategories of $D^b(\proj\,R_i)$,
for all $i$. Using this information, how can we get a classification of all triangulated subcategories of $D^b(\proj\,R)$?  
More generally, one can ask how the spectral theory of $R$ and the spectral theories of the rings $R_i$ are related.

Before we go further, we remark that by the obvious induction, it suffices to consider only two components 
($R \cong R_1 \times  R_2)$. We collect some standard facts abouts products of triangulated categories
that will be needed.

\begin{lemma} \label{products} Let $\T_1$ and $\T_2$ be triangulated categories. Then we have the following.
\begin{itemize}
\item The product category $\T_1 \times \T_2$ admits a triangulated structure that has the following 
universal property: Given any triangulated category $\D$ and some triangulated functors 
$\T_1 \leftarrow \D \rightarrow \T_2$, there exist a unique triangulated functor $\F: \D \rightarrow \T_1 \times \T_2$ 
making the following diagram of triangulated functors commutative. 
\[
\xymatrix{
& \D \ar[dl] \ar@{.>}[d]^F \ar[dr] &\\
\T_1 & \T_1 \times \T_2 \ar[l]^{\pi_1} \ar[r]_{\pi_2} & \T_2
}
\]

 \item Every thick (localising) subcategory of $\T_1 \times \T_2$ is of the form $\B_1 \times \B_2$, where
$\B_i$ is a thick (localising) subcategory in $\T_i$.

\item If both $\T_1$ and $\T_2$ are essentially small triangulated 
(tensor triangulated) categories, then $K_0(\T_1 \times \T_2) \cong K_0(\T_1) \times K_0(\T_2)$ as groups (rings).
\end{itemize}
\end{lemma}

\begin{rem} In contrast with the thick subcategories, not every triangulated subcategory of
$\T_1 \times \T_2$ is of the form $\A_1 \times \A_2$, where $\A_i$ is a triangulated subcategory 
of $\T_i$. This is clear from the third part of the above lemma because not every subgroup of
the product group $K_0(\T_1) \times K_0(\T_2)$ is a product of subgroups. 
\end{rem}

Now we relate the category of perfect complexes over $R$ and those over the rings
$R_i$.

\begin{prop} There is a natural equivalence of triangulated categories,
\[ D^b(\proj\,R) \simeq D^b(\proj\,R_1) \times D^b(\proj\,R_2).\]
In particular $K_0[D^b(\proj\,R)]  \cong K_0 [D^b(\proj\,R_1)] \times K_0[D^b(\proj\,R_2)].$
\end{prop}

\begin{proof} Note that every module $M$ over $R_1 \times R_2$ is a direct sum (in the category of
$R_1 \times R_2$ - modules) 
\[ M \cong P_1 \oplus P_2,\]
where $P_i$ is an $R_i$-module: Take $P_1 = \la (1,0) \ra M$ and $P_2 = \la (0,1) \ra M$, where $P_i$ is also regarded as 
an $R_1 \times R_2$ module via the projection maps $R_1 \times R_2 \rightarrow R_i$. Now one can verify easily that 
this decomposition is functorial and sends finitely generated (projective) modules to finitely generated (projective) 
modules. Therefore every complex in $D^b(\proj\,R_1 \times R_2)$ splits as a direct sum of two complexes, one each in 
$D^b(\proj\,R_i)$. Conversely, given a pair of complexes $(X_1, X_2)$ with $X_i$ in $D^b(\proj\,R_i)$, their direct sum 
$X_1 \oplus X_2$ is clearly a complex in $D^b(\proj\, R_1 \times R_2)$. Now it can be verified  that these two functors 
establish the desired equivalence of triangulated categories.  The second statement follows immediately from the 
lemma \ref{products}.
\end{proof}

So in view of this proposition, the problem of classifying triangulated subcategories in $D^b(\proj\,R_1 \times R_2)$
boils down to classifying triangulated subcategories of $D^b(\proj\,R_1) \times D^b(\proj\,R_2)$. So, following
Thomason's recipe, we need to classify the thick subcategories of $D^b(\proj\,R_1) \times D^b(\proj\,R_2)$ and compute their Grothendieck groups. 
This is given by Lemma \ref{products}: Every thick subcategory $\T$ of $D^b(\proj\,R_1) \times D^b(\proj\,R_2)$
is of the form $\T_1 \times \T_2$ where $\T_i$ is thick in $D^b(\proj\,R_i)$, and 
$K_0(\T) \cong K_0(\T_1 \times \T_2) \cong K_0(\T_1) \times K_0(\T_2)$.

\subsubsection{Artin rings}
\noindent
We will now apply these ideas to Artin rings.  Recall that a ring is said to be \emph{Artinian} if every descending chain of ideals terminates. 
Artin rings can be characterised as zero dimensional noetherian rings. They have the following structure theorem. 
\begin{thm}\cite{am}
Every Artin ring $R$ is isomorphic to a finite direct product of 
      Artin local rings. Moreover, the number of local rings that appear in this
      isomorphism is equal to the cardinality of $\Spec(R)$.
\end{thm}

Thus  $R \cong \prod_{i=1}^{n} R_i$, where each $R_i$ is an Artin local ring
($n$ is the cardinality of $\Spec(R)$). We have seen above that there is an equivalence 
of triangulated categories,
\[ D^b(\proj\,R) \cong \prod_{i=1}^{n} D^b(\proj\,R_i).\]
So by the above discussion, we just have to compute the Grothendieck groups of the thick subcategories of $D^b(\proj\,R_i)$.
But since each $R_i$ is an Artinian local ring, $\Spec(R_i)$ is a one-point space. This implies (by the Hopkins-Neeman theorem) that
the only non-zero thick subcategory is $D^b(\proj\,R_i)$ itself, whose Grothendieck group is well-known to be infinite cyclic.  
Therefore we have,
 \[ K_0(\T_S) \cong \bigoplus_{p_i \in S} K_0(D^b(\proj\,R_i)) \cong \bigoplus_{p_i \in S} \mathbb{Z}. \]
The universal Euler characteristic function that gives this isomorphism is $\oplus_{p_i \in S} \Lambda_{p_i}$ where 
$\Lambda_{p_i} (X) = \sum_{t=-\infty}^{\infty} (-1)^t \dim_{R_i/p_i} H_t (X \otimes R_i/p_i) $. 

We now record the classification of triangulated subcategories of perfect complexes over Artin rings.

\begin{thm} Let $R$ be any Artin ring and let $R = \prod_i R_i$ be its unique decomposition into Artin local rings. 
For every subset $S$ of $\Spec(R)$ and every subgroup  $H$ of $\bigoplus_{p_i \in S} \mathbb{Z}$, define
\[\T(S,H) := \{X \in \T_S: (\oplus_{p_i \in S} \Lambda_{p_i} )(X) \in H \}. \] 
This is a complete list of  triangulated subcategories of $D^b(\proj \, R)$.  
Further, every dense triangulated subcategory of $\T_S$ is a triangulated ideal if and only if $S$ is a
one point space.
\end{thm}

\begin{rem} It is clear from the proof that this theorem also holds whenever $R \cong \prod_{i=1}^{n} R_i$, where each of the rings 
$R_i$ has exactly one prime ideal.  
\end{rem}

We now derive some easy consequences of the above theorem.

\begin{cor}\label{cor:questionBArtin}
Let $X$ and $Y$ be perfect complexes over an Artin ring. Then $Y$ can be generated from
$X$ using cofibrations if and only if
\begin{itemize}
\item $\Supp(Y) \subseteq \Supp(X)$, and 
\item $\Lambda_{p_i}(X)$ divides $\Lambda_{p_i}(Y)$ for all $p_i \in \Supp(X)$,
\end{itemize}
\end{cor}

\begin{cor} An Artin ring $R$ is local if and only if every dense triangulated subcategory of $D^b(\proj\,R)$ is  a triangulated ideal.
\end{cor}

\begin{proof} By theorem \ref{imain}, it is clear that every dense triangulated subcategory is a triangulated ideal if and only
if every subgroup of $K_0(D^b(\proj\,R)) \cong \prod_{p \in \Spec(R)} \ints$ is also an ideal. Clearly the later happens if and only if 
$|\Spec(R)| = 1$, or equivalently if $R$ is local.
\end{proof}

\noindent
\subsection{Non-noetherian rings}

In order to study the problem of classifying  triangulated subcategories in the non-noetherian case,
we need a thick subcategory theorem for $D^b(\proj\,R)$, when $R$ is a non-noetherian ring.
This is given by a result of Thomason, which is a far-reaching generalisation of 
Hopkins-Neeman theorem to schemes. We now state this theorem for commutative rings.

\begin{thm} \cite{Th} \label{cor:thicksubcatrings} Let $R$ be any commutative ring. Then there is a natural order preserving bijection between
the sets
\begin{center}
\{thick subcategories $\A$ of $D^b(\proj \, R)$\}
\begin{center}
$f\downarrow  \;\;\; \uparrow g$
\end{center}
\{subsets $S$ of $\Spec(R)$ such that $ S = \bigcup_{\alpha} V(I_{\alpha})$, where $I_{\alpha}$ is finitely generated\}.
\end{center}
The map $f$ sends a thick subcategory $\A$ to the set $ \bigcup_{X \in \A} \Supp(X) $
and the map $g$ sends $S$ to the thick subcategory $\{ X \in D^b(\proj \, R): \Supp(X) \in S \}$. 
\end{thm}

\begin{rem} The subsets of $\Spec(R)$ in this corollary which determine the thick subcategories of perfect complexes will be
called \emph{Thick supports}.  If $R$ is noetherian, every thick support is a specialisation-closed subset (since ideals in a noetherian ring
are finitely generated), therefore the above corollary recovers the Hopkins-Neeman thick subcategory theorem.
\end{rem}

We do not know much about the $K$-theory of thick subcategories over non-noetherian rings, except in the simplest case where
the rings have only one prime ideal, e.g., $R = \ftwo[X_2, X_3, \cdots ]/(X_2^2, X_3^3, \cdots )$.

\begin{rem}
The geometry of these rings is very simple -- just a one point space, however, their derived categories can be  very mysterious
and extremely complicated.  Amnon Neeman \cite{ne2} showed that the derived category of the above ring has uncountably many Bousfield classes --
a striking contrast with the noetherian result where the Bousfield classes are known to be in bijection with the subsets
of $\Spec(R)$ \cite{Ne}. Despite this incredible complexity in the derived categories of such rings, $K$-theory does enable us to classify all
the triangulated subcategories of perfect complexes. 
\end{rem}

\begin{prop} Let $R$ be any commutative ring with a unique prime ideal $p$. Then every triangulated subcategory of
$D^b(\proj\,R)$ is a triangulated ideal and is of the form
  \[ \D_m = \{X \in D^b(\proj\,R)\; \mbox{such that}\; \Lambda(X) \equiv 0 \!\!\! \mod{m} \} \]
for some non-negative integer $m$, where 
$\Lambda(X) = \sum_{-\infty}^\infty (-1)^i \dim_{R/p} H_i(X \otimes R/p)$.
\end{prop}
\begin{proof} Note that any such $R$ is, in particular, a local ring. Therefore, 
\[K_0(D^b(\proj\,R)) (\cong K_0(R)) \cong \mathbb{Z}.\]
It is easily verified that the given Euler characteristic function gives this isomorphism. 
Moreover, since $R$ has a unique prime ideal, it is clear from corollary \ref{cor:thicksubcatrings} that there are no non-trivial
thick subcategories in $D^b(\proj\,R)$. So the dense triangulated
subcategories in $D^b(\proj\,R)$ are all the triangulated subcategories in $D^b(\proj\,R)$. It is clear that these are all triangulated ideals.
This completes the proof of the proposition.
\end{proof}

\subsection{Questions}

\subsubsection{Algebraic $K$-theory for thick subcategories.}
The key to the problem of classifying the triangulated subcategories of perfect complexes lies in the algebraic $K$-theory of thick subcategories.
In section \ref{se:algktheory} we have summarised a few results from classical algebraic $K$-theory. Those
were results about the Grothendieck groups of $D^b(\proj\,R)$. So we now ask if such results also hold for thick subcategories
of $D^b(\proj\,R)$. We ask a very specific question to make this point clear. 
It is well-known that if $J$ is the nilradical of $R$, then
\[ K_0(D^b(\proj\,R)) \cong K_0(D^b(\proj\,R/J)).\]
Since every prime ideal contains the nilradical, the quotient map $R \rar R/J$ induces a homeomorphism on prime spectra:
$\Spec(R) \cong \Spec(R/J)$. This homeomorphism of prime spectra implies that the lattice of specialisation-closed subsets of 
$\Spec(R)$ is isomorphic to that of $\Spec(R/J)$. Now if $R$ is noetherian, we can invoke the Hopkins-Neeman thick subcategory theorem to conclude that 
the same is true for the lattices of thick subcategories of perfect complexes over $R$ and $R/J$. Now the question that arises is whether
the thick subcategories that correspond to each other under this isomorphism have isomorphic Grothendieck groups?


\end{document}

%% file: ah.pstex_t
\begin{picture}(0,0)%
\includegraphics{ah.pstex}%
\end{picture}%
\setlength{\unitlength}{1618sp}%
\begingroup\makeatletter\ifx\SetFigFont\undefined%
\gdef\SetFigFont#1#2#3#4#5{%
  \reset@font\fontsize{#1}{#2pt}%
  \fontfamily{#3}\fontseries{#4}\fontshape{#5}%
  \selectfont}%
\fi\endgroup%
\begin{picture}(16847,10241)(-2999,-6373)
\end{picture}%

%% file: bss.pstex_t
\begin{picture}(0,0)%
\includegraphics{bss.pstex}%
\end{picture}%
\setlength{\unitlength}{3315sp}%
\begingroup\makeatletter\ifx\SetFigFont\undefined%
\gdef\SetFigFont#1#2#3#4#5{%
  \reset@font\fontsize{#1}{#2pt}%
  \fontfamily{#3}\fontseries{#4}\fontshape{#5}%
  \selectfont}%
\fi\endgroup%
\begin{picture}(8424,10224)(2989,-9973)
\end{picture}%

%% file: chromatic.bbl
\begin{thebibliography}{BCR97}

\bibitem[Ada69]{adams}
J.~F. Adams.
\newblock Lectures on generalised cohomology.
\newblock In {\em Category Theory, Homology Theory and their Applications, III
  (Battelle Institute Conference, Seattle, Wash., 1968, Vol. Three)}, pages
  1--138. Springer, Berlin, 1969.

\bibitem[Ada92]{MayTho}
J.~F. Adams.
\newblock {\em The selected works of {J}. {F}rank {A}dams. {V}ol. {II}}.
\newblock Cambridge University Press, Cambridge, 1992.

\bibitem[AK89]{adamskuhn}
J.~F. Adams and N.~J. Kuhn.
\newblock Atomic spaces and spectra.
\newblock {\em Proc. Edinburgh Math. Soc. (2)}, 32(3):473--481, 1989.

\bibitem[AM69]{am}
M.~F. Atiyah and I.~G. Macdonald.
\newblock {\em Introduction to commutative algebra}.
\newblock Addison-Wesley Publishing Co., Reading, Mass.-London-Don Mills, Ont.,
  1969.

\bibitem[BCR97]{bcr}
D.~J. Benson, Jon~F. Carlson, and Jeremy Rickard.
\newblock Thick subcategories of the stable module category.
\newblock {\em Fund. Math.}, 153(1):59--80, 1997.

\bibitem[Boa99]{boardman}
J.~Michael Boardman.
\newblock Conditionally convergent spectral sequences.
\newblock In {\em Homotopy invariant algebraic structures (Baltimore, MD,
  1998)}, volume 239 of {\em Contemp. Math.}, pages 49--84. Amer. Math. Soc.,
  Providence, RI, 1999.

\bibitem[Chr98]{ch}
J.~Daniel Christensen.
\newblock Ideals in triangulated categories: phantoms, ghosts and skeleta.
\newblock {\em Adv. Math.}, 136(2):284--339, 1998.

\bibitem[DHS88]{dhs}
Ethan~S. Devinatz, Michael~J. Hopkins, and Jeffrey~H. Smith.
\newblock Nilpotence and stable homotopy theory. {I}.
\newblock {\em Ann. of Math. (2)}, 128(2):207--241, 1988.

\bibitem[Fau03]{fau}
H.~Fausk.
\newblock Picard groups of derived categories.
\newblock {\em J. Pure Appl. Algebra}, 180(3):251--261, 2003.

\bibitem[Gro77]{sga5}
A.~Grothendieck.
\newblock {\em Cohomologie {$l$}-adique et fonctions {$L$}}.
\newblock Springer-Verlag, Berlin, 1977.

\bibitem[Hop87]{Ho}
Michael~J. Hopkins.
\newblock Global methods in homotopy theory.
\newblock In {\em Homotopy theory (Durham, 1985)}, volume 117 of {\em London
  Math. Soc. Lecture Note Ser.}, pages 73--96. Cambridge Univ. Press,
  Cambridge, 1987.

\bibitem[Hov99]{hoveymodel}
Mark Hovey.
\newblock {\em Model categories}, volume~63 of {\em Mathematical Surveys and
  Monographs}.
\newblock American Mathematical Society, Providence, RI, 1999.

\bibitem[HPS97]{mps}
Mark Hovey, John~H. Palmieri, and Neil~P. Strickland.
\newblock Axiomatic stable homotopy theory.
\newblock {\em Mem. Amer. Math. Soc.}, 128(610):x+114, 1997.

\bibitem[HS98]{hs}
Michael~J. Hopkins and Jeffrey~H. Smith.
\newblock Nilpotence and stable homotopy theory. {II}.
\newblock {\em Ann. of Math. (2)}, 148(1):1--49, 1998.

\bibitem[Mit85]{mit}
Stephen~A. Mitchell.
\newblock Finite complexes with {$A(n)$}-free cohomology.
\newblock {\em Topology}, 24(2):227--246, 1985.

\bibitem[Nee92]{Ne}
Amnon Neeman.
\newblock The chromatic tower for {$D(R)$}.
\newblock {\em Topology}, 31(3):519--532, 1992.

\bibitem[Nee00]{ne2}
Amnon Neeman.
\newblock Oddball {B}ousfield classes.
\newblock {\em Topology}, 39(5):931--935, 2000.

\bibitem[Rav84]{rav}
Douglas~C. Ravenel.
\newblock Localization with respect to certain periodic homology theories.
\newblock {\em Amer. J. Math.}, 106(2):351--414, 1984.

\bibitem[Tho97]{Th}
R.~W. Thomason.
\newblock The classification of triangulated subcategories.
\newblock {\em Compositio Math.}, 105(1):1--27, 1997.

\bibitem[Wei94]{wei}
Charles~A. Weibel.
\newblock An introduction to homological algebra.
\newblock 38:xiv+450, 1994.

\bibitem[Wei03]{wei1}
Charles~A. Weibel.
\newblock {\em Algebraic K-theory}.
\newblock http://math.rutgers.edu/\~{}weibel, 2003.

\bibitem[Xu95]{xu}
Kai Xu.
\newblock A note on atomic spectra.
\newblock In {\em The \v Cech centennial (Boston, MA, 1993)}, volume 181 of
  {\em Contemp. Math.}, pages 419--422. Amer. Math. Soc., Providence, RI, 1995.

\end{thebibliography}
